\documentstyle{amsppt}
\pagewidth{15.5truecm} \pageheight{20truecm}
\tolerance =100000
\NoRunningHeads
\TagsOnRight
\headline ={\tenrm\hfil\folio\hfil}

\def\sqr#1#2{{\vcenter{\vbox{\hrule height.#2pt
              \hbox{\vrule width.#2pt height#1pt \kern#1pt \vrule width.#2pt}
              \hrule height.#2pt}}}}
\def\signed #1{{\unskip\nobreak\hfil\penalty50
              \hskip2em\hbox{}\nobreak\hfil#1
              \parfillskip=0pt \finalhyphendemerits=0 \par}}
\def\endpf{\signed {$\sqr69$}}

\def\no{\noindent}

\def\med{\medskip}

\def\IE{\text{\rm l\negthinspace E}}

\def\IN{\text{\rm l\negthinspace N}}

\def\IR{\text{\rm l\negthinspace R}}

\topmatter
\title\nofrills
Dynamic programming principle for one kind of stochastic recursive
optimal control problem  and Hamilton-Jacobi-Bellman equations*
\thanks * This work is supported by
the National Natural Science Foundation (10671112) and the New
century  Young Teachers Program of Education Ministry , P.R.China.
\endthanks
\endtitle
\author Zhen WU$^a$ \thanks $^a$ Email: wuzhen\@sdu.edu.cn\endthanks
\quad  zhiyong YU$^b$ \thanks $^b$ Email:
yuzhiyong\@sdu.edu.cn\endthanks
\endauthor
\affil School of Mathematics and System Sciences\\
 Shandong University\\
  Jinan 250100\hskip2.00mm  P.R.China
\endaffil
\endtopmatter
\document
\no{\bf Abstract.} In this paper, we study one kind of stochastic
recursive optimal control problem with the obstacle constraints
for the cost function where the cost function is  described by the
solution of one reflected backward stochastic differential
equations.  We will give the dynamic programming principle for
this kind of optimal control problem and show that the value
function is the unique viscosity solution of the obstacle problem
for the corresponding Hamilton-Jacobi-Bellman equations.
\par\

\no {\bf Keywords:} Reflected backward stochastic differential
equation, Recursive optimal control problem, Dynamic programming
principle, Hamilton-Jacobi-Bellman equations, Viscosity solution.
\par\

\no {\bf AMS subject classification:} 93E20, 60H10, 35K15 \par\

\no{\bf 1. Introduction.} \med

Nonlinear backward stochastic differential equations (BSDE in short)
have
  been  introduced by Pardoux \& Peng [11].
Independently,
  Duffie \& Epstein
  [6] introduced BSDE from economic background. In [6] they presented a
stochastic  differential recursive utility which is an extension of
the standard additive utility with the instantaneous utility
depending not only on the instantaneous consumption rate $c_t$ but
also on the future utility. Actually it corresponds to the solution
of a particular BSDE associated with a generator  which does not
depend on the variable $z$. In mathematics the result in [11] is
more general.  Then El.Karoui, Peng and Quenez [10] gave some
important properties of BSDE such as comparison theorem and
applications in mathematical finance and optimal control theory. And
also in this paper they gave the formulation of recursive utilities
and their properties from the BSDE point of view. The recursive
optimal control problem is presented as a kind of optimal control
problem whose cost function is described by the solution of BSDE. In
1992, Peng [12]  got the Bellman's dynamic programming principle for
this kind of problem and proved the value function is a viscosity
solution of one kind of quasi-linear second-order partial
differential equation (PDE in short) which is the well-known
Hamilton-Jacobi-Bellman equation. Then in 1997, he virtually
generalized these results to much more general situation, even under
Non-Markvian framework. (See in [13]). In this chinese version, Peng
used the backward semigroup property of BSDE to give a complete
proof of the Bellman's dynamic programming principle for the
recursive optimal problem introduced by a BSDE whose coefficient
just satisfies Lipschitz condition, under Markovian and
Non-Markovian framework.  He also proved that the value function is
a viscosity solution of a generalized Hamilton-Jacobi-Bellman
equation. \med

  Then El.Karoui, Kapoudjian, Pardoux, Peng and Quenez [9] studied the
   reflected BSDE with one barrier. The solution of the reflected BSDE is forced
   to stay above one given continuous stochastic process which is called
   ``obstacle". For this purpose they introduced one increasing process to push
 the solution upwards in a kind of minimal way. They got
 the existence and uniqueness of the solution for this kind of reflected BSDE
  and also studied its
 relation with the obstacle problem for nonlinear parabolic PDE's within the Markov
 framework. Using two different methods, Snell envelope theory connected with
 fixed point principle and penalization method. Cvitanic and Karatzas [5]
 extended the result to reflected BSDE's with two barriers
 called upper and lower barriers, which are two given
 continuous processes.
   Hamad\`ene and Lepeltier [7] generalized the results of El.Karoui
et al [9] to one barrier which is right continuous and left upper
semi-continuous. They used this model to solve the mixed optimal
stochastic control problem when the terminal reward is only right
continuous and  left upper semi-continuous. In this kind of mixed
control problem, the controller has two
 actions,
one is of control and the other is of stopping his control strategy
in view to maximize his payoff. Also in this paper Hamad\`ene and
Lepeltier generalized the result of Cvitanic and Karatzas ([5]) to
reflected BSDE's with two barriers to processes $S$ (lower barrier)
and $-U$ ($U$ is upper barrier)  merely right continuous and left
upper semicontinuous.  And then Hamad\`ene, Lepeltier and Wu [8]
proved existence and uniqueness results of the solution for infinite
horizon
 reflected backward stochastic differential equations with one or two barriers.
 They  also apply those results to get the existence of optimal control strategy
 for the mixed control problem and a saddle-point strategy for the mixed game
 problem when, in both situation, the  horizon is infinite.

\med

   In our paper, we study one kind of recursive optimal control
   problem with the  obstacle constraints for the cost function. This means that the
   cost function of the control system is described by the
   solution of one reflected BSDE which is required to satisfy the
   obstacle constraints. This kind of the recursive optimal
   control problem has some practical sense such as, in financial
   market, the investor requires his recursive utility function
   value to be bigger than one specific function of his wealth. For this
   purpose, one increasing process is introduced to push the cost
   function value upward and we also hope this push power to be
   minimum. From the result in [9] and [7], we know that, in
   fact, this kind of problem is one mixed recursive optimal
   stochastic control problem.\med

     One of our interesting problem is that if the dynamic
     programming principle still holds for the above optimal control
     problem. Using some properties of the reflected BSDE and
     analysis technique we give the positive answer for this
     question. This result can be seen as the generalized
     extension of the dynamic programming principle of the
     recursive control problem in [12] and [13] to the obstacle
     constraints case for the cost function. And then, we show that, provided the
     problem is formulated within a Markovian framework, the value
     function is the unique viscosity solution of the obstacle
     problem for one nonlinear parabolic PDEs which is called Hamilton-Jacobi-Bellman
     (HJB in short) equations.\med
     The paper is organized as follows. In section 2, we present
     some preliminary results about reflected stochastic
     differential equations which play important role to study the
     dynamic programming principle of the optimal control problem.
     In section 3, we formulate the recursive optimal control
     problem with the obstacle constraints for the cost function and prove that the dynamic programming
     principle still holds. In section 4, we show that the value function of
     the control problem is the unique viscosity solution of the
     obstacle problem for corresponding HJB equations. In Appendix we
     put in some  technique proof of the preliminary results of the reflected BSDE.

\med \no{\bf 2. Preliminary results of the reflected BSDE} \med

 In this section, we give some preliminary results of the
 reflected BSDE which is useful to get the dynamic programming
 principle for the recursive optimal control problem with the
 obstacle constraints for the cost function.

 Let $\{W_t, 0\leq t\leq T\}$ be a $d-$dimensional standard
Brownian motion defined on a probability space $(\Omega, \Cal F,
P)$. Let $\{\Cal F_t, 0\leq t\leq T\}$ be the natural filtration
of $\{W_t\}$, where $\Cal F_0$ contains all P-null sets of $\Cal
F$ and let $\Cal P$ be the $\sigma-$algebra of predictable subsets
of $\Omega\times [0,T]$.

Let us introduce some notation.
$$\aligned
 L^2 &=
\left\{\text{$\xi$ is an $\Cal F_T-$ measurable random variable
s.t.
} \IE(|\xi|^2)<+\infty\right\},\\
H^2 &= \left\{ \{\varphi_t, 0\leq t\leq T \} \text { is a predictable process s.t. } \IE\int_0^T|\varphi_t|^2dt<+\infty \right\},\\
S^2 &= \left\{ \{\varphi_t, 0\leq t\leq T \} \text{ is a
predictable process s.t. }\IE(\sup_{0\leq t\leq
T}|\varphi_t|^2)<+\infty \right\}\endaligned$$ and the following
reflected BSDE with one barrier:
$$ Y_t=\xi+\int_t^T g(s,Y_s,Z_s)ds+K_T-K_t-\int_t^T Z_sdW_s,\qquad 0\leq t\leq T. \tag 2.1$$
Here $\xi\in L^2$, $g$ is a map from $\Omega\times [0,T]\times
\IR\times\IR^d$ onto $\IR$ satisfying

(i) $\forall$ $(y,z)\in\IR\times\IR^d$, $g(\cdot,y,z)\in H^2$,

(ii) for some $L>0$ and all $y,y'\in\IR$, $z,z'\in\IR^d$, a.s.
$$|g(t,y,z)-g(t,y',z')|\leq L(|y-y'|-|z-z'|),$$
an ``obstacle" $\{S_t, 0\leq t\leq T\}$, which is a continuous
progressively measurable real-valued process satisfying

(iii) $\IE\left(\sup_{0\leq t\leq T}|S_t|^2\right)<+\infty$.

Then from Theorem 5.2 in [9], there exists unique solution
$\{(Y_t,Z_t,K_t), 0\leq t\leq T\}$  taking values in $\IR$, $\IR^d$
and $\IR_+$, respectively, and satisfying:

(iv) $Y\in S^2$, $Z\in H^2$ and $K_T\in L^2$;

(v) $Y_t\geq S_t,\qquad 0\leq t\leq T;$

(vi) $\{K_t\}$ is continuous and increasing, $K_0=0$ and
$$\int_0^T(Y_t-S_t)dK_t=0.$$

Now we  give two more accurate estimates on the norm of the solution
similar to Proposition 3.5 and Proposition 3.6 in [9].\med

\no{\bf Proposition 2.1} \quad {\it Let $\{(Y_t,Z_t,K_t), 0\leq t
\leq T\}$ be the solution of the above reflected BSDE, then there
exists a constant $C$ such that}

$$
\IE^{\Cal F_t}\left\{\sup_{t\leq s\leq T }
Y_s^2+\int_t^T|Z_s|^2+|K_T-K_t|^2\right\}\leq C\IE^{\Cal
F_t}\left\{\xi^2+\left(\int_t^T g(s,0,0)ds\right)^2+\sup_{t\leq
s\leq T} S_t^2\right\}. $$

This proposition is similar to Proposition 3.5 in [9]. However, the
estimate is more precise which is necessary to get the desired
results in next section. The proof is a little complicated and
technical, some technique derive from [2], we put it in the
Appendix. \med

And then, we need to  estimate the variation of  the solution
induced by a variation of the reflected BSDE coefficients.\med

\no{\bf Proposition 2.2} \quad {\it Let $(\xi,g,S)$ and
$(\xi^{\prime},g^{\prime},S^{\prime})$ be two triplets satisfying
the above assumptions. Suppose $(Y,Z,K)$ is the solution of the
reflected BSDE $(\xi,g,S)$ and
$(Y^{\prime},Z^{\prime},K^{\prime})$ is the solution of the
reflected BSDE $(\xi^{\prime},g^{\prime},S^{\prime})$. Define
$$\Delta\xi=\xi-\xi^{\prime},\qquad \Delta g=g-g^{\prime},\qquad \Delta S=S-S^{\prime};$$
$$\Delta Y=Y-Y^{\prime},\qquad \Delta Z=Z-Z^{\prime},\qquad \Delta K=K-K^{\prime}.$$
Then there exists a constant $C$ such that}
$$\aligned
&\IE^{\Cal F_t}\left\{\sup_{t\leq s\leq T}|\Delta Y_s|^2
+\int_t^T|\Delta Z_s|^2ds +|\Delta K_T-\Delta K_t|^2\right\}\\
&\leq C\IE^{\Cal F_t}\left\{|\Delta\xi|^2+\left(\int_t^T |\Delta
g(s,Y_s,Z_s)|ds\right)^2\right\}+C\left(\IE^{\Cal
F_t}\left\{\sup_{t\leq s\leq T}|\Delta
S_s|^2\right\}\right)^{1/2}\Psi_{t,T}^{1/2},
\endaligned$$
{\it where} $$\aligned \Psi_{t,T} &=\IE^{\Cal
F_t}\left\{|\xi|^2+\left(\int_t^T|g(s,0,0)|ds
\right)^2 +\sup_{t\leq s\leq T} |S_s|^2\right.\\
&\quad \left.+|\xi^{\prime}|^2+\left(\int_t^T|g^{\prime}(s,0,0)|ds
\right)^2 +\sup_{t\leq s\leq T} |S^{\prime}_s|^2\right\}.
\endaligned$$

 The estimate of this proposition is more accurate than that in
 Proposition 3.6 in [9]. We also put the proof in the
 Appendix.
\med \no{\bf 3. Formulation of the problem and Dynamic programming
principle} \med

In this section, we first formulate one kind of stochastic
recursive optimal control problem with the obstacle constraints
for the cost function, and then we prove that dynamic programming
principle still holds for this kind of optimization problem.
\med

   We introduce the admissible control set $\Cal U$ defined by
$$\Cal U :=\left\{ v(\cdot)\in H^2 |\quad v(\cdot)\quad \text{take value in } U\subset\IR^k\right\}.$$
$U$ is a compact set, the element of $\Cal U$ is called admissible
control.

 For given admissible control, we consider the following control system
$$\left\{\aligned
 dX^{t,\zeta;v}_s &=
b(s,X^{t,\zeta;v}_s,v_s)ds+\sigma(s,X^{t,\zeta;v}_s,v_s)dW_s,\qquad
s\in [t,T],\\
X^{t,\zeta;v}_t &= \zeta,
\endaligned\right.\tag 3.1$$
here $t\ge 0$ is regarded as the initial time, $\zeta\in
L^2(\Omega,\Cal F_t,P;\IR^n)$ as the initial state, the mappings
$$b : [0,T]\times\IR^n\times U\to \IR^n,\qquad \sigma : [0,T]\times\IR^n\times U\to \IR^{n\times d}$$
satisfy the following conditions:

(H3.1) $b$ and $\sigma$ are continuous in $t$;

(H3.2) For some $L>0$, and all $x,x^{\prime}\in\IR^n$,
$v,v^{\prime}\in U$, a.s.
$$\aligned
&|b(t,x,v)-b(t,x^{\prime},v^{\prime})|+|\sigma(t,x,v)-\sigma(t,x^{\prime},v^{\prime})|\leq
L(|x-x^{\prime}|+|v-v^{\prime}|).\endaligned$$

Obviously, under above assumptions, for any $v(\cdot)\in \Cal U$,
control system (3.1) has a unique strong solution
$\{X^{t,\zeta;v}_s$, $0\le t\le s\le T\}$,  and we also have the
following estimates:\med

\no{\bf Proposition 3.1} \quad {\it For all $t\in [0,T]$,
$\zeta,\zeta^{\prime}\in L^2(\Omega,\Cal F_t,P;\IR^n)$,
$v(\cdot),v^{\prime}(\cdot)\in\Cal U$,
$$\IE^{\Cal F_t}\left\{\sup_{t\leq s\leq T}|X^{t,\zeta;v}_s|^2\right\}\leq C(1+|\zeta|^2);\tag 3.2$$
$$\IE^{\Cal F_t}\left\{\sup_{t\leq s\leq T}|X^{t,\zeta;v}_s-X^{t,\zeta';v'}_s|^2\right\}
\leq C|\zeta-\zeta^{\prime}|^2+C\IE^{\Cal
F_t}\left\{\int_t^T|v_s-v^{\prime}_s|^2ds\right\},\tag 3.3$$ where
the constant $C$ depends only on $L$.}\med

\no{\bf Proposition 3.2} \quad {\it For all $t\in [0,T]$, $x\in
\IR^n$, $v(\cdot)\in \Cal U$, $\delta\in[0,T-t]$,
$$\IE\left\{\sup_{t\leq s\leq t+\delta}|X^{t,x;v}_s-x|^2\right\}\leq C\delta,\tag 3.4$$
where the constant $C$ depend only on $x$ and $L$}.\med

Now  for any given admissible control $v(\cdot)\in\Cal U$, we
consider the following reflected BSDE
$$\aligned Y^{t,\zeta;v}_s&=\Phi(X^{t,\zeta;v}_T)+\int_s^Tg(r,X^{t,\zeta;v}_r,Y^{t,\zeta;v}_r,Z^{t,\zeta;v}_r,v_r)dr\\
&\qquad
+K^{t,\zeta;v}_T-K^{t,\zeta;v}_s-\int_s^TZ^{t,\zeta;v}_rdW_r,\qquad
t\leq s\leq T,\endaligned\tag 3.5$$ here
$$\aligned
& \Phi=\Phi(x) : \IR^n\to\IR,\quad h=h(t,x) :
[0,T]\times\IR^n\to\IR,\\
& g=g(t,x,y,z,v) : [0,T]\times \IR^n\times\IR\times\IR^d\times
U\to \IR\endaligned$$
 satisfy the following conditions:

(H3.3)  $g$ and $h$ are continuous in $t$;

(H3.4) For some $L>0$, and all $x,x^{\prime}\in\IR^n$,
$y,y^{\prime}\in\IR$, $z,z^{\prime}\in\IR^d$ $v,v^{\prime}\in U$,
a.s. $$\aligned
&|g(t,x,y,z,v)-g(t,x^{\prime},y^{\prime},z^{\prime},v^{\prime})|+|\Phi(x)-\Phi(x^{\prime})|+|h(t,x)-h(t,x^{\prime})|\\
&\leq
L(|x-x^{\prime}|+|y-y^{\prime}|+|z-z^{\prime}|+|v-v^{\prime}|).\endaligned$$
Then from Theorem 5.2 in [9],
 there exists a
unique triple $(Y^{t,\zeta;v},Z^{t,\zeta;v},K^{t,\zeta;v})$, which
is the solution of reflected BSDE (3.5), satisfying

(i) $Y^{t,\zeta;v}\in S^2$, $Z^{t,\zeta;v}\in H^2$  and
$K^{t,\zeta;v}_T\in L^2$;

(ii) $Y^{t,\zeta;v}_s\geq h(s,X^{t,\zeta;v}_s)$, $t\leq s\leq T$;

(iii) $\{K^{t,\zeta;v}_s\}$  is increasing and continuous,
$K^{t,\zeta;v}_t=0$,  and $
\int_t^T(Y^{t,\zeta;v}_s-h(s,X^{t,\zeta;v}_s))dK^{t,\zeta;v}_s=0$.

Moreover, we can get the following estimates for the solution of
(3.5) from  Proposition 2.1 and 2.2.\med

\no{\bf Proposition 3.3}
$$\IE^{\Cal F_t}\left\{\sup_{t\leq s\leq T}|Y^{t,\zeta;v}_s|^2+\int_t^T|Z^{t,\zeta;v}_s|^2ds+|K^{t,\zeta;v}_T|^2\right\}\leq C(1+|\zeta|^2).\tag 3.6$$
\med \no{\bf Proposition 3.4}
$$\aligned
&\IE^{\Cal F_t}\left\{\sup_{t\leq s\leq
T}|Y^{t,\zeta;v}_s-Y^{t,\zeta^{\prime};v^{\prime}}|^2+\int_t^T|Z^{t,\zeta;v}_s-Z^{t,\zeta^{\prime};v^{\prime}}_s|^2ds+
|K^{t,\zeta;v}_T-K^{t,\zeta^{\prime};v^{\prime}}_T|^2\right\}\\
&\quad\leq C|\zeta-\zeta^{\prime}|^2+C\IE^{\Cal
F_t}\left\{\int_t^T|v_s-v^{\prime}_s|^2ds\right\}\\
&\qquad+C(1+|\zeta|+|\zeta^{\prime}|)\left(|\zeta-\zeta^{\prime}|^2+\IE^{\Cal
F_t}\left\{\int_t^T|v_s-v^{\prime}_s|^2ds\right\}\right)^{1/2}.\endaligned\tag
3.7$$ \med

Given the control process $v(\cdot)\in\Cal U$, we introduce the
associated cost functional:
$$J(t,x;v(\cdot)):= Y^{t,x;v}_s|_{s=t},\qquad (t,x)\in [0,T]\times\IR^n,\tag 3.8$$
and we define the value function of the stochastic optimal control
problem
$$u(t,x):=ess\sup_{v(\cdot)\in\Cal U} J(t,x;v(\cdot)),\qquad (t,x)\in [0,T]\times\IR^n.\tag 3.9$$
\med

 This is one kind of stochastic recursive optimal control problem
 with the obstacle constraints for the cost function:
$Y^{t,x;v}_s\ge h(s,X^{t,x;v}_s)$, $t\le s\le T$. In financial
market, if $X^{t,x;v}_s$ represents the wealth of the one investor
, $Y^{t,x;v}_s$: the recursive utility cost function, the
constraint is that the investor requires his cost function value
to be bigger than one function of his wealth at any time.\med

\no{\bf Remark 3.5} $\quad$  From Proposition 2.3 in [9] and the
definition in [7] and [8], we know that the above optimal control
problem is one recursive mixed optimal control problem:
$$\aligned
u(t,x):=ess\sup_{v(\cdot)\in\Cal U}
Y^{t,x;v}_t&=ess\sup_{v(\cdot)\in\Cal U}ess\sup_{\tau\in\Cal
{T}_t}\IE^{\Cal F_t}\left\{\int_t^{\tau}g(s,X^{t,\zeta;v}_s,Y^{t,\zeta;v}_s,Z^{t,\zeta;v}_s,v_s)ds\right.\\
&\quad
\left.+h(X^{t,\zeta;v}_{\tau})1_{\tau<T}+\Phi(X^{t,\zeta;v}_T)1_{\tau=T}\right\}\endaligned$$
where $\Cal T$ is the set of all stopping times dominated by $T$
and $\Cal T_t=\{\tau\in\Cal T;\quad t\le\tau\le T\}$.

  In this kind of recursive mixed control problem, the controller has two
  actions, one is of control $v(\cdot)$ and the other is of
  stopping his control strategy in view to maximize his recursive
  payoff. The more detail about this kind of problem can be seen
  in [9], [7] and [8].
  \med

  Now we continue to study the former control problem (3.9) and
  show that celebrated dynamic programming principle still holds
  for this kind of optimization problem. The main proof idea  comes from  the proof of
  dynamic programming principle for recursive problem given by Peng in chinese version [13].
 \med

For each $t>0$, we denote by $\{\Cal F^t_s, t\leq s\leq T\}$ the
natural filtration of the Brownian motion $\{W_s-W_t, t\leq s\leq
T\}$, augmented  by the P-null sets of $\Cal F$ and we introduce
the following subspaces of admissible controls
$$\aligned
&\Cal U^t := \left\{ v(\cdot)\in\Cal U\quad |\quad v(s)\text{ is
$\{\Cal F^t_s\}$ progressively measurable, } \forall\ t\leq s\leq
T.
\right\}\\
&\bar{\Cal U}^t := \left\{ v_s=\sum_{j=1}^N v^j_s1_{A_j}\quad
|\quad v^j_s\in\Cal U^t,\quad \{A_j\}_{j=1}^N\ \text {is a
partition of } (\Omega,\Cal F_t). \right\}
\endaligned$$
\med Firstly we will show that
\med
 \no{\bf Proposition 3.6}
{\it Under the assumptions (H3.1)--(H3.4), the value function
$u(t,x)$ defined in (3.9) is a deterministic function.} \med {\bf
Proof}:  Firstly , we will show
$$ess\sup_{v(\cdot)\in\Cal U}J(t,x;v(\cdot))=ess\sup_{v(\cdot)\in\bar{\Cal U}^t}J(t,x;v(\cdot)).\tag 3.10$$
Obviously,
$$ess\sup_{v(\cdot)\in\Cal U}J(t,x;v(\cdot))\geq ess\sup_{v(\cdot)\in\bar{\Cal U}^t}J(t,x;v(\cdot)).$$
we need to show the inverse inequality. $\forall\varepsilon>0$,
there exists $\tilde{v}(\cdot)\in \Cal U$ such that
$$P\left\{J(t,x;\tilde{v}(\cdot))>ess\sup_{v(\cdot)\in\Cal U}J(t,x;v(\cdot))-\varepsilon\right\}=\delta>0.$$
From (3.7), we know $\forall\bar{v}(\cdot)\in\bar{\Cal U}^t$,
$$\IE\left\{|Y^{t,x;\bar{v}}_t-Y^{t,x;\tilde{v}}_t|^2\right\}\leq C\IE\int_t^T|\bar{v}_s-\tilde{v}_s|^2ds+C(1+x)\left
(\IE\int_t^T|\bar{v}_s-\tilde{v}_s|^2ds\right)^{1/2}.$$
 Note that
$\bar{\Cal U}^t$ is dense in $\Cal U$, then there exists a
sequence $\{v_n(\cdot)\}_{n=1}^\infty\in\bar{\Cal U}^t$ such that
$$\lim_{n\to\infty}\IE\left\{|Y^{t,x;v_n}_t-Y^{t,x;\tilde{v}}_t|^2\right\}=0.$$
Then, there exists a subsequence, we denote without loss of
generality $\{v_n(\cdot)\}_{n=1}^\infty$ also, such that
$$\lim_{n\to\infty}Y^{t,x;v_n}_t=Y^{t,x;\tilde{v}}_t\quad a.s.,$$
then
$$\aligned
& P\left(\bigcap_{m=1}^\infty\bigcup_{N=1}^\infty\bigcap_{n=N}^\infty\left\{|Y^{t,x;v_n}_t-Y^{t,x;\tilde{v}}_t|<\frac{1}{m}
\right\}\right)=1,\\
& P\left(\bigcup_{N=1}^\infty\bigcap_{n=N}^\infty\left\{|Y^{t,x;v_n}_t-Y^{t,x;\tilde{v}}_t|<\frac{1}{m}\right\}\right)=1,
\qquad \forall m\in\IN,\\
& \lim_{N\to\infty}P\left(\bigcap_{n=N}^\infty\left\{|Y^{t,x;v_n}_t-Y^{t,x;\tilde{v}}_t|<\frac{1}{m}\right\}\right)=1,
\qquad \forall m\in\IN,\\
&
\lim_{N\to\infty}P\left\{|Y^{t,x;v_n}_t-Y^{t,x;\tilde{v}}_t|<\frac{1}{m}\right\}=1,\qquad
\forall m\in\IN.
\endaligned$$
We select $m$ big enough such that $1/m<\varepsilon$ and denote
$$\aligned
A &= \left\{\omega | Y^{t,x;\tilde{v}}_t>ess\sup_{v(\cdot)\in\Cal
U}J(t,x;v(\cdot))-\varepsilon\right\};\\
B_N &= \left\{\omega |
|Y^{t,x;v_N}_t-Y^{t,x;\tilde{v}}_t|\leq\frac{1}{m}\right\},\quad
N=1,2,\cdots,
\endaligned$$
then, from above definition, $P(A)=\delta>0$ and
$\lim_{N\to\infty}P(B_N)=1$. We select $N$ big enough such that
$P(B_N)>1-\delta$, then
$$P(AB_N)=P(A)+P(B_N)-P(A\cup B_N)>\delta+(1-\delta)-1=0.$$
It is easily to check
$$P\left\{Y^{t,x;v_N}_t>ess\sup_{v(\cdot)\in\Cal U}J(t,x;v(\cdot))-2\varepsilon\right\}\geq P(AB_N)>0.$$
This inequality implies
$$ess\sup_{v(\cdot)\in\bar{\Cal U}^t}J(t,x;v(\cdot))\geq ess\sup_{v(\cdot)\in\Cal U}J(t,x;v(\cdot))-2\varepsilon.$$
From the arbitrariness of $\varepsilon$, we get
$$ess\sup_{v(\cdot)\in\bar{\Cal U}^t}J(t,x;v(\cdot))\geq ess\sup_{v(\cdot)\in\Cal U}J(t,x;v(\cdot)).$$
Then we obtain (3.10).

Secondly, we will show
$$ess\sup_{v(\cdot)\in\bar{\Cal U}^t}J(t,x;v(\cdot))=ess\sup_{v(\cdot)\in\Cal U^t}J(t,x;v(\cdot))\tag 3.11$$
Obviously,
$$ess\sup_{v(\cdot)\in\bar{\Cal U}^t}J(t,x;v(\cdot))\geq ess\sup_{v(\cdot)\in\Cal U^t}J(t,x;v(\cdot)).$$
We need to show the inverse inequality also.

Let us admit for a moment the following lemma. The main idea of the
lemma is to consider the partition of probability space, which is
first introduced by Theorem 4.7 in [13].   \med

\no{\bf Lemma 3.7}
$$\aligned
X^{t,x;\sum_{j=1}^Nv^j1_{A_j}}_. &=
\sum_{j=1}^N1_{A_j}X^{t,x;v^j}_.;\quad
Y^{t,x;\sum_{j=1}^Nv^j1_{A_j}}_. = \sum_{j=1}^N1_{A_j}Y^{t,x;v^j}_.;\\
Z^{t,x;\sum_{j=1}^Nv^j1_{A_j}}_. &=
\sum_{j=1}^N1_{A_j}Z^{t,x;v^j}_.;\quad
K^{t,x;\sum_{j=1}^Nv^j1_{A_j}}_. =
\sum_{j=1}^N1_{A_j}K^{t,x;v^j}_..
\endaligned$$
\med \no $\forall v(\cdot)\in\bar{\Cal U}^t$, we have
$$J(t,x;v(\cdot))=J(t,x;\sum_{j=1}^Nv^j(\cdot)1_{A_j})=\sum_{j=1}^N1_{A_j}J(t,x;v^j(\cdot)).$$
Note that $v^j(\cdot)$ are $\{\Cal F^t_s\}$ progressively
measurable, then $J(t,x;v^j(\cdot))$ $(j=1,2,\cdots,N)$ are
deterministic. Without loss of generality, we assume that
$$J(t,x;v^1(\cdot))\geq J(t,x;v^j(\cdot)),\quad \forall j=2,3,\cdots,N.$$
So that
$$J(t,x;v(\cdot))\leq J(t,x;v^1(\cdot))\leq ess\sup_{v(\cdot)\in\Cal U^t}J(t,x;v(\cdot)).$$
From the arbitrariness of $v(\cdot)$, we get
$$ess\sup_{v(\cdot)\in\bar{\Cal U}^t}J(t,x;v(\cdot))\leq ess\sup_{v(\cdot)\in\Cal U^t}J(t,x;v(\cdot)),$$
and obtain (3.11).

However, when $v(\cdot)\in\Cal U^t$, the cost functional
$J(t,x;v(\cdot))$ is deterministic, so
$$u(t,x)=\sup_{v(\cdot)\in\Cal U^t}J(t,x;v(\cdot))$$
is deterministic and the  proof is completed.
\endpf
\med

We need to give
\med
\no{\bf Proof of Lemma 3.7} $\quad$   For every $j=1,2,\cdots,N$,
we denote
$$(X^j_s, Y^j_s, Z^j_s, K^j_s) \equiv (X^{t,x;v^j}_s, Y^{t,x;v^j}_s, Z^{t,x;v^j}_s, K^{t,x;v^j}_s).$$
$X^j$ is the solution of the following stochastic differential
equations:
$$X^j_s=x+\int_t^sb(r,X^j_r,v^j_r)dr+\int_t^s\sigma(r,X^j_r,v^j_r),\quad s\in [t,T].$$
$(Y^j,Z^j,K^j)$ satisfies the following reflected BSDE:
$$\aligned
& Y^j_s=\Phi(X^j_T)+\int_s^Tg(r,X^j_r,Y^j_r,Z^j_r,v_r)dr+K^j_T-K^j_s-\int_s^TZ^j_rdW_r,\quad s\in [t,T];\\
& Y^j_s\geq h(s,X^j_s),\quad s\in [t,T];\quad
\int_t^T(Y^j_s-h(s,X^j_s))dK^j_s=0.
\endaligned$$
We multiply $1_{A_j}$ on the both sides of the above equations,
then sum the equations. From the  trivial fact:
$$\sum_j1_{A_j}\varphi(x_j)=\varphi(\sum_jx_j1_{A_j}),$$
we get
$$\aligned
\sum_{j=1}^N1_{A_j}X^j_s&=x+\int_t^sb(r,\sum_{j=1}^N1_{A_j}X^j_r,\sum_{j=1}^N1_{A_j}v^j_r)dr
+\int_t^s\sigma(r,\sum_{j=1}^N1_{A_j}X^j_r,\sum_{j=1}^N1_{A_j}v^j_r)dW_r;\\
\sum_{j=1}^N1_{A_j}Y^j_s&=\Phi(\sum_{j=1}^N1_{A_j}X^j_T)+\int_s^Tg(r,\sum_{j=1}^N1_{A_j}X^j_r,\sum_{j=1}^N1_{A_j}Y^j_r,\sum_{j=1}^N1_{A_j}Z^j_r,\sum_{j=1}^N1_{A_j}v^j_r)dr\\
&\quad+\sum_{j=1}^N1_{A_j}K^j_T-\sum_{j=1}^N1_{A_j}K^j_s-\int_s^T\sum_{j=1}^N1_{A_j}Z^j_r;\\
\sum_{j=1}^N1_{A_j}Y^j_s&\geq
h(s,\sum_{j=1}^N1_{A_j}X^j_s);\quad\int_t^T\left(\sum_{j=1}^N1_{A_j}Y^j_s-h(s,\sum_{j=1}^N1_{A_j}X^j_s)\right)d\left(\sum_{j=1}^N1_{A_j}K^j_s\right)=0.
\endaligned$$
Then from the uniqueness of the solution of stochastic
differential equations  and reflected BSDE, we get the conclusion.
\endpf
\med We next will discuss the continuity of value function
$u(t,x)$ with respect to $x$. We have the following
estimation:\med

\no{\bf Lemma 3.8}  {\it For each $t\in [0,T]$, $x$ and
$x^{\prime}\in\IR^n$, we have}

 (i) $|u(t,x)-u(t,x^{\prime})|^2\leq C|x-x^{\prime}|^2+C(1+|x|+|x^{\prime}|)|x-x^{\prime}|;$

 (ii)$|u(t,x)|\leq C(1+|x|).$

{\bf Proof}: $\quad$  From estimation (3.6) and (3.7), for each
admissible control $v(\cdot)\in\Cal U$, we have
$$|J(t,x;v(\cdot))|\leq C(1+|x|); $$
$$ |J(t,x;v(\cdot))-J(t,x^{\prime};v(\cdot))|^2\leq
C|x-x^{\prime}|^2+C(1+|x|+|x^{\prime}|)|x-x^{\prime}|.\tag 3.12
$$
On the other hand, for each $\varepsilon>0$, there exist
$v(\cdot)$ and $v^{\prime}(\cdot)\in\Cal U$ such that
$$\aligned
J(t,x;v^{\prime}(\cdot))&\leq u(t,x)\leq
J(t,x;v(\cdot))+\varepsilon,\\
J(t,x';v(\cdot))&\leq u(t,x')\leq
J(t,x';v'(\cdot))+\varepsilon.\endaligned$$ Then from the
estimation of $J$, we get
$$-C(1+|x|)\leq J(t,x;v^{\prime}(\cdot))\leq u(t,x)\leq J(t,x;v(\cdot))+\varepsilon\leq C(1+|x|)+\varepsilon.$$
From the arbitrariness of $\varepsilon$, we obtain (ii).
Similarly,
$$J(t,x;v^{\prime}(\cdot))-J(t,x^{\prime};v^{\prime}(\cdot))-\varepsilon\leq u(t,x)-u(t,x')\leq J(t,x;v(\cdot))
-J(t,x^{\prime};v(\cdot))+\varepsilon,$$
$$\aligned
& |u(t,x)-u(t,x^{\prime})| \\
&\leq
\max\left\{|J(t,x;v(\cdot))-J(t,x^{\prime};v(\cdot))|,|J(t,x;v^{\prime}(\cdot))
-J(t,x^{\prime};v^{\prime}(\cdot))|\right\}+\varepsilon,\endaligned$$
$$\aligned
& |u(t,x)-u(t,x^{\prime})|^2\\
&\leq
2\max\left\{|J(t,x;v(\cdot))-J(t,x^{\prime};v(\cdot))|^2,|J(t,x;v'(\cdot))-J(t,x^{\prime};
v^{\prime}(\cdot))|^2\right\}+2\varepsilon^2\\
&\leq
2C|x-x^{\prime}|^2+2C(1+|x|+|x^{\prime}|)|x-x^{\prime}|+2\varepsilon^2.
\endaligned$$
The we can obtain (i).
\endpf
\med

We also have\med

\no{\bf Lemma 3.9} $\quad$ {\it $\forall t\in [0,T]$, $\forall
v(\cdot)\in\Cal U$, for all $\zeta\in L^2(\Omega,\Cal
F_t,P;\IR^n)$, we have
$$J(t,\zeta;v(\cdot))=Y^{t,\zeta;v}_t.$$
}

{\bf Proof}:  We first study the simple case: $\zeta$ has the
following form:
$$\zeta=\sum_{i=1}^N 1_{A_i}x_i,$$
where $\{A\}_{i=1}^N$ is a finite partition of $(\Omega,\Cal
F_t)$, and $x_i\in \IR^n$, for $1\leq i\leq N$. The similar
argument as Lemma 3.7 leads to
$$Y^{t,\zeta;v}_s=Y^{t,\sum_{i=1}^N1_{A_i}x_i;v}_s=\sum_{i=1}^N1_{A_i}Y^{t,x_i;v}_s,\qquad s\in [t,T].$$
From the definition (3.8), we deduce that
$$Y^{t,\zeta;v}_t=\sum_{i=1}^N1_{A_i}Y^{t,x_i;v}_t=\sum_{i=1}^N1_{A_i}J(t,x_i;v(\cdot))
=J(t,\sum_{i=1}^N1_{A_i}x_i;v(\cdot))=J(t,\zeta;v(\cdot)).$$
Therefore, for simple functions, we have the desired result.

Given a general $\zeta\in L^2(\Omega,\Cal F_t,P;\IR^n)$, we can
choose a sequence of simple functions $\{\zeta_i\}$ which
converges to $\zeta$ in $L^2(\Omega,\Cal F_t,P;\IR^n)$.
Consequently, from the estimate (3.7) and (3.12), we have
$$\aligned
&\IE\left\{|Y^{t,\zeta;v}_t-Y^{t,\zeta_i;v}_t|^2\right\}\\
&\leq
\IE\left\{C|\zeta-\zeta_i|^2+C(1+|\zeta|+|\zeta_i|)|\zeta-\zeta_i|\right\}\\
&\leq
C\IE\left\{|\zeta-\zeta_i|^2\right\}+C\left(\IE\left\{(1+|\zeta|+|\zeta_i|)^2\right\}\right)^{1/2}\left(\IE\left\{
|\zeta-\zeta_i|^2\right\}\right)^{1/2}\\
&\rightarrow 0,\qquad as\ i\to\infty,\endaligned$$
$$\aligned
&\IE\left\{|J(t,\zeta;v(\cdot))-J(t,\zeta_i;v(\cdot))|^2\right\}\\
&\leq
\IE\left\{C|\zeta-\zeta_i|^2+C(1+|\zeta|+|\zeta_i|)|\zeta-\zeta_i|\right\}\\
&\leq
C\IE\left\{|\zeta-\zeta_i|^2\right\}+C\left(\IE\left\{(1+|\zeta|+|\zeta_i|)^2\right\}\right)^{1/2}\left(\IE\left\{
|\zeta-\zeta_i|^2\right\}\right)^{1/2}\\
&\rightarrow 0,\qquad as\ i\to\infty,
\endaligned$$
and $Y^{t,\zeta_i;v}_t=J(t,\zeta_i;v(\cdot))$, the proof is
completed.
\endpf
\med

  For the value function of our recursive optimal control problem,
  we have
\med \no{\bf Lemma 3.10} $\quad$ {\it Fixed $t\in[0,T)$ and
$\zeta\in L^2(\Omega,\Cal{F}_t,P;\IR^n)$, for each $v(\cdot)\in
\Cal{U}$, we have
$$
u(t,\zeta)\geq Y^{t,\zeta;v}_t.\tag 3.13
$$
On the other hand, for each $\varepsilon>0$, there exists an
admissible control $v(\cdot)\in\Cal{U}$ such that
$$
u(t,\zeta)\leq Y^{t,\zeta;v}_t+\varepsilon, \quad a.s..\tag 3.14
$$
}
\med
{\bf Proof}:  We first prove (3.13). When $\zeta$ is a simple
function:
$$\zeta=\sum_{i=1}^N1_{A_i}x_i,$$
for all $v(\cdot)\in \Cal U$, we have
$$Y^{t,\zeta;v}_t=Y^{t,\sum_{i=1}^N1_{A_i}x_i;v}_t=\sum_{i=1}^N1_{A_i}Y^{t,x_i;v}_t\leq
\sum_{i=1}^N1_{A_i}u(t,x_i)=u(t,\zeta). $$

When $\zeta\in L^2(\Omega,\Cal F_t,P;\IR^n)$, we can choose a
sequence of simple functions $\{\zeta_i\}$ which converges to
$\zeta$ in $L^2(\Omega,\Cal F_t,P;\IR^n)$. Consequently, similarly
with Lemma 3.9, we have
$$
\IE\left\{|Y^{t,\zeta;v}_t-Y^{t,\zeta_i;v}_t|^2\right\}\rightarrow
0;\quad\IE\left\{|u(t,\zeta)-u(t,\zeta_i)|^2\right\}\rightarrow 0.
$$ Then, there exists a subsequence, we use same notation without
loss of generality also, such that
$$
\lim_{i\to\infty}Y^{t,\zeta_i;v}_t=Y^{t,\zeta;v}_t,\qquad
a.s.\quad\lim_{i\to\infty}u(t,\zeta_i)=u(t,\zeta),\qquad a.s. $$
here  $Y^{t,\zeta_i;v}_t\leq u(t,\zeta_i)$, $i=1,2,\cdots$, so
$Y^{t,\zeta;v}_t\leq u(t,\zeta)$.

We turn to prove (3.14). We first deal with the case that $\zeta$
is a bounded random variable: $\zeta\in L^\infty(\Omega,\Cal
F_t,P;\IR^n)$. We suppose that $|\zeta|\leq M$ and construct a
simple random variable $\eta\in L^\infty(\Omega,\Cal F_t,P;\IR^n)$
$$\eta=\sum_{i=1}^N1_{A_i}x_i$$
such that

(i) $|\eta|\leq |\zeta|$;

 (ii) $|\eta-\zeta|\leq
\min\left\{\frac{\varepsilon}{6\sqrt{C}},\frac{\varepsilon^2}{36C(1+2M)}\right\}$.

 For any $v(\cdot)\in\Cal U$, we have
$$
|Y^{t,\zeta;v}_t-Y^{t,\eta;v}_t| \leq \frac{\varepsilon}{3};\quad
|u(t,\zeta)-u(t,\eta)| \leq \frac{\varepsilon}{3}. $$ Then for
each $x_i$, we can choose an $\{\Cal F^t_s\}-$adapted admissible
control $v^i(\cdot)$ such that
$$u(t,x_i)\leq Y^{t,x_i;v_i}+\frac{\varepsilon}{3}.$$
We denote
$$v(\cdot):=\sum_{i=1}^N1_{A_i}v^i(\cdot),$$
then
$$\aligned
 Y^{t,\zeta;v}_t &\geq
-|Y^{t,\zeta;v}_t-Y^{t,\eta;v}_t|+Y^{t,\eta;v}_t \geq
-\frac{\varepsilon}{3}+\sum_{i=1}^N1_{A_i}Y^{t,x_i;v_i}_t\\
&\geq -\frac \varepsilon 3
+\sum_{i=1}^N1_{A_i}(u(t,x_i)-\frac{\varepsilon}{3})
= -\frac 2 3 \varepsilon +u(t,\eta)\\
&\geq -\varepsilon+u(t,\zeta).
\endaligned$$
Therefore, for $\zeta\in L^\infty(\Omega,\Cal F_t,P;\IR^n)$, we
have the desired result (3.14).

Given a general $\zeta\in L^2(\Omega,\Cal F_t,P;\IR^n)$, we note
that $\zeta$ have the following form:
$$\zeta=\sum_{i=1}^\infty1_{A_i}\zeta_i,$$
where $\{A_i\}_{i=1}^\infty$ is a partition of $(\Omega,\Cal
F_t)$, $x_i\in\IR^n$ $(i=1,2,\cdots)$, $|\zeta_i|\leq i$ and
$\zeta_i\in L^\infty(\Omega,\Cal F_t,P;\IR^n)$. So, for every
$\zeta_i$, there exists $v^i(\cdot)\in\Cal U$, such that
$$u(t,\zeta_i)\leq Y^{t,\zeta_i;v_i}+\varepsilon.$$
We denote
$$v(\cdot)=\sum_{i=1}^\infty1_{A_i}v^i(\cdot),$$
and then
$$\aligned
u(t,\zeta) &= u(t,\sum_{i=1}^\infty1_{A_i}\zeta_i) =
\sum_{i=1}^\infty1_{A_i}u(t,\zeta_i) \leq
\sum_{i=1}^\infty1_{A_i}(Y^{t,\zeta_i;v_i}+\varepsilon)\\
&=\sum_{i=1}^\infty1_{A_i}Y^{t,\zeta_i;v_i}+\varepsilon =
Y^{t,\zeta;v}_t+\varepsilon.
\endaligned$$
The proof is completed.
\endpf
\med

Now we start to discuss the (generalized) dynamic programming
principle for our recursive optimal control problem (3.9). In [13],
Peng first used the idea of (backward) semigroups of BSDE to prove
the dynamic programming principle for the recursive optimal control
problem associated to BSDE.\med

Firstly we introduce a family of (backward) semigroups which come
from Peng's idea [13].

Given the initial condition $(t,x)$, an admissible control
$v(\cdot)\in\Cal U$, a positive number $\delta\leq T-t$ and a
real-valued random variable $\eta\in L^2(\Omega,\Cal
F_{t+\delta},P;\IR)$, we denote
$$G^{t,x;v}_{t,t+\delta}[\eta]:=Y_t,$$
where $(Y_s,Z_s,K_s)_{t\leq s\leq t+\delta}$ is the solution of
the following reflected BSDE with time horizon $t+\delta$
$$\aligned
Y_s&=\eta+\int_s^{t+\delta}g(r,X^{t,x;v}_r,Y_r,Z_r,v_r)dr+K_{t+\delta}-K_s\\
& \qquad -\int_s^{t+\delta}Z_rdW_r,\qquad t\leq s\leq
t+\delta,\endaligned$$ satisfying

(i) $Y\in S^2$, $Z\in H^2$  and $K_{t+\delta}\in L^2$;

(ii) $Y_s\geq h(s,X^{t,x;v}_s)$, $t\leq s\leq t+\delta$;

(iii) $\{K_s\}$  is increasing and continuous,  $K_t=0$,
$\int_t^{t+\delta}(Y_s-h(s,X^{t,x;v}_s))dK_s=0$.

 Obviously,
$$G^{t,x;v}_{t,T}[\Phi(X^{t,x;v}_T)]=G^{t,x;v}_{t,t+\delta}[Y^{t,x;v}_{t+\delta}].$$

Then our (generalized) dynamic programming principle holds.\med

\no{\bf Theorem 3.11}  {\it Under the assumptions (H3.1)--(H3.4),
the value function $u(t,x)$ obeys the following dynamic
programming principle: For each $0<\delta\leq T-t$,
$$u(t,x)=ess\sup_{v(\cdot)\in\Cal U} G^{t,x;v}_{t,t+\delta}[u(t+\delta,X^{t,x;v}_{t+\delta})].\tag 3.15$$
}

{\bf Proof}:  We have
$$\aligned
u(t,x) &= ess\sup_{v(\cdot)\in\Cal
U}G^{t,x;v}_{t,T}[\Phi(X^{t,x;v}_T)]= ess\sup_{v(\cdot)\in\Cal
U}G^{t,x;v}_{t,t+\delta}[Y^{t,x;v}_{t+\delta}]\\
&= ess\sup_{v(\cdot)\in\Cal
U}G^{t,x;v}_{t,t+\delta}[Y^{t+\delta,X^{t,x;v}_{t+\delta};v}_{t+\delta}].
\endaligned$$
From Lemma 3.10 and the comparison theorem of reflected BSDE
(Theorem 4.1 in [9]),
$$u(t,x)\leq  ess\sup_{v(\cdot)\in\Cal
U}G^{t,x;v}_{t,t+\delta}[u(t+\delta,X^{t,x;v}_{t+\delta})].$$ On
the other hand, for every $\varepsilon>0$, we can find an
admissible control $\bar{v}(\cdot)\in\Cal U$ such that
$$u(t+\delta,X^{t,x;v}_{t+\delta})\leq Y^{t+\delta,X^{t,x;v}_{t+\delta};\bar{v}}_{t+\delta}+\varepsilon.$$
From this and the comparison theorem, we get
$$u(t,x)\geq ess\sup_{v(\cdot)\in\Cal U}G^{t,x;v}_{t,t+\delta}[u(t+\delta,X^{t,x;v}_{t+\delta})-\varepsilon].$$
From Proposition 2.2, there exists a positive constant $C_0$ such
that
$$u(t,x)\geq ess\sup_{v(\cdot)\in\Cal U}G^{t,x;v}_{t,t+\delta}[u(t+\delta,X^{t,x;v}_{t+\delta})]-C_0\varepsilon.$$
Therefore, letting $\varepsilon\downarrow 0$, we obtain the
equation (3.15).
\endpf
\med At the end of this section, we devote ourselves to obtaining
the continuity of $u(t,x)$ with respect to $t$.\med

\no{\bf Proposition 3.12}  {\it The value function $u(t,x)$ is
continuous in t.}
\med
{\bf Proof}:  We define $Y^{t,x;v}_s$ for all $s\in [0,T]$ by
choosing $Y^{t,x;v}_s\equiv Y^{t,x;v}_t$ for $0\leq s\leq t$. And
we define the ``obstacle"
$$S^{t,x;v}_s=\left\{\aligned &h(s,X^{t,x;v}_s) ; t\leq s\leq T;\\ &h(t,x) ; \qquad\quad 0\leq s\leq t. \endaligned
\right.$$

Fixed $x\in\IR^n$, for all $0\leq t_1\leq t_2\leq T$, we analysis
the difference of $u(t_1,x)$ and $u(t_2,x)$ .

$\forall \varepsilon>0$, there exist $v_1(\cdot)\in\Cal U$,
$v_2(\cdot)\in\Cal U$, such that
$$
 Y^{t_1,x;v_2}_{t_1}\leq u(t_1,x)\leq
Y^{t_1,x;v_1}_{t_1}+\varepsilon;\quad Y^{t_2,x;v_1}_{t_2}\leq
u(t_2,x)\leq Y^{t_2,x;v_2}_{t_2}+\varepsilon. $$ Then,
$$Y^{t_1,x;v_2}_{t_1}-Y^{t_2,x;v_2}_{t_2}-\varepsilon\leq u(t_1,x)-u(t_2,x)\leq Y^{t_1,x;v_1}_{t_1}-Y^{t_2,x;v_1}_{t_2}+\varepsilon,$$
$$|u(t_1,x)-u(t_2,x)|\leq\max\{|Y^{t_1,x;v_1}_{t_1}-Y^{t_2,x;v_1}_{t_2}|,
|Y^{t_1,x;v_2}_{t_1}-Y^{t_2,x;v_2}_{t_2}|\}+\varepsilon.$$ Here we
only estimate $|Y^{t_1,x;v_1}_{t_1}-Y^{t_2,x;v_1}_{t_2}|$ and the
estimate of $|Y^{t_1,x;v_2}_{t_1}-Y^{t_2,x;v_2}_{t_2}|$ is same.
From Proposition 2.2, we have
$$\aligned
&|Y^{t_1,x;v_1}_{t_1}-Y^{t_2,x;v_1}_{t_2}|^2=|Y^{t_1,x;v_1}_0-Y^{t_2,x;v_1}_0|^2\\
&\leq \IE\left\{\sup_{0\leq s\leq T}|Y^{t_1,x;v_1}_s-Y^{t_2,x;v_1}_s|^2\right\}\\
&\leq C\IE\left\{|\Phi(X^{t_1,x;v_1}_T)-\Phi(X^{t_2,x;v_1}_T)|^2\right\}\\
&\quad +C\IE\left\{\left(\int_0^T|1_{[t_1,
T]}g(s,X^{t_1,x;v_1}_s,Y^{t_1,x;v_1}_s,Z^{t_1,x;v_1}_s,v_1(s))\right.\right.\\
& \quad\qquad\qquad -1_{[t_2,
T]}g(s,X^{t_2,x;v_1}_s,Y^{t_1,x;v_1}_s,Z^{t_1,x;v_1}_s,v_1(s))|ds\big)^2\big\}\\
&\quad +C\Psi^{1/2}_{0,T}\left(\IE\left\{\sup_{0\leq s\leq T}
|S^{t_1,x;v_1}_s-S^{t_2,x;v_1}_s|^2\right\}\right)^{1/2},
\endaligned\tag 3.16$$

where
$$\aligned \Psi_{0,T} &=
\IE\left\{|\Phi(X^{t_1,x;v_1}_T)|^2+\left(\int_{t_1}^T|g(s,X^{t_1,x;v_1}_s,0,0,v_1(s))|ds\right)^2\right.\\
& \quad +\sup_{t_1\leq s\leq
T}|h(s,X^{t_1,x;v_1}_s)|^2+|\Phi(X^{t_2,x;v_1}_T)|^2\\
& \quad
\left.+\left(\int_{t_2}^T|g(s,X^{t_2,x;v_1}_s,0,0,v_1(s))|ds\right)^2+\sup_{t_2\leq
s\leq T}|h(s,X^{t_2,x;v_1}_s)|^2\right\}.
\endaligned$$
Now we deal with the items for the right side of inequality
(3.16).

The first item: From Lipschitz condition, Proposition 3.1 and
Proposition 3.2, we get
$$I\leq C\IE\left\{|X^{t_1,x;v_1}_T-X^{t_2,x;v_1}_T|^2\right\}\leq C\IE\{|X^{t_2,x;v_1}_{t_1}-x|^2\}\leq C(t_2-t_1).$$

The second item: From Lipschitz condition, $(a+b)^2\leq
a^2/2+b^2/2$,  Proposition 3.1, Proposition 3.2 and Proposition
3.3, we get
$$II\leq C(t_2-t_1).$$

The third item: As the same argument we get
$$\Psi_{0,T}\leq C.$$
We next discuss
$$\aligned
|S^{t_1,x;v_1}_s-S^{t_2,x;v_1}_s|^2&=|h(s,X^{t_1,x;v_1}_s)-h(s,X^{t_2,x;v_1}_s)|^2\\
&\leq C|X^{t_1,x;v_1}_s-X^{t_2,x;v_1}_s|^2;\qquad s\in
[t_2,T],\endaligned$$
$$\aligned
|S^{t_1,x;v_1}_s-S^{t_2,x;v_1}_s|^2&=|h(s,X^{t_1,x;v_1}_s)-h(t_2,x)|^2\\
&\leq C|X^{t_1,x;v_1}_s-x|^2+2|h(s,x)-h(t_2,x)|^2;\qquad s\in
[t_1,t_2],\endaligned$$
$$
|S^{t_1,x;v_1}_s-S^{t_2,x;v_1}_s|^2=|h(t_1,x)-h(t_2,x)|^2;\qquad
s\in [0,t_1].$$

So we have $$\aligned &\IE\left\{\sup_{0\leq s\leq
T}|S^{t_1,x;v_1}_s-S^{t_2,x;v_1}_s|^2\right\}\\
&\leq \IE\left\{\left(\sup_{0\leq s\leq t_1}+\sup_{t_1\leq s\leq
t_2}+\sup_{t_2\leq s\leq
T}\right)|S^{t_1,x;v_1}_s-S^{t_2,x;v_1}_s|^2\right\}\\
&\leq C(t_2-t_1)+|h(t_1,x)-h(t_2,x)|^2+2\sup_{t_1\leq s\leq
t_2}|h(s,x)-h(t_2,x)|^2\\
&\leq C(t_2-t_1)+3\sup_{t_1\leq s\leq t_2}|h(s,x)-h(t_2,x)|^2.
\endaligned$$
From the above analysis , we know
$$|Y^{t_1,x;v_1}_{t_1}-Y^{t_2,x;v_1}_{t_2}|^2\leq C(t_2-t_1)+3\sup_{t_1\leq s\leq t_2}|h(s,x)-h(t_2,x)|^2,$$
$$|Y^{t_1,x;v_1}_{t_1}-Y^{t_2,x;v_1}_{t_2}|\leq C(t_2-t_1)^{1/2}+3\sup_{t_1\leq s\leq t_2}|h(s,x)-h(t_2,x)|.$$
The same argument used to
$|Y^{t_1,x;v_2}_{t_1}-Y^{t_2,x;v_2}_{t_2}|^2$ leads to
$$|u(t_1,x)-u(t_2,x)|\leq C(t_2-t_1)^{1/2}+3\sup_{t_1\leq s\leq t_2}|h(s,x)-h(t_2,x)|+\varepsilon.$$
Because of the arbitrariness of $\varepsilon$, we get
$$|u(t_1,x)-u(t_2,x)|\leq C(t_2-t_1)^{1/2}+3\sup_{t_1\leq s\leq t_2}|h(s,x)-h(t_2,x)|.$$
From the continuity of $h(t,x)$ with respect to $t$, we get the
continuity of $u(t,x)$ with respect to $t$. The proof is
completed.
\endpf
\med

\no{\bf 4. Viscosity solution of an obstacle problem for HJB
equations} \med

 In this section, we relate the value function of above recursive
 optimal control problem with the following obstacle problem for
 nonlinear second-order parabolic PDEs
 which is called Hamilton-Jacobi-Bellman equations:

$$\left\{
\aligned
& \min\big(u(t,x)-h(t,x),\\
&\quad\quad -\frac{\partial u}{\partial t}(t,x)
-\sup_{v\in U}\left\{\Cal L(t,x,v)u(t,x)+g(t,x,u(t,x),\nabla u(t,x)\sigma(t,x,v),v)\right\}\big)=0,\\
& u(T,x)=\Phi(x),
\endaligned\right.\tag 4.1$$
where $\Cal L$ is a family of second order linear partial
differential operators,
$$\Cal L(t,x,v)\varphi=\frac{1}{2}Tr\left((\sigma\sigma^T)(t,x,v)D^2\varphi\right)+\langle b(t,x,v), D\varphi\rangle.$$
Here the function $b,\sigma,g,\Phi,h$ are supposed to satisfy
(H3.1)--(H3.4), respectively. \med We want to prove that the value
function $u(t,x)$ introduced by (3.9) is the unique viscosity
solution of the obstacle problem for HJB equation (4.1). We first
recall the definition of a viscosity solution for HJB equation
obstacle problem (4.1) from [4]. Below, $S^n$ will denote the set of
$n\times n$ symmetric matrices.\med

\no{\bf Definition 4.1} {\it Let $u(t,x)\in C((0,T)\times\IR^n)$
and $(t,x)\in (0,T)\times\IR^n$. We denote by $\Cal P^{2,+}u(t,x)$
[the ``parabolic superjet" of $u$ at $(t,x)$] the set of triples
$(p,q,X)\in\IR\times\IR^n\times S^n$ which are such that
$$
u(s,y) \leq u(t,x)+p(s-t)+\langle q, y-x\rangle
        +\frac{1}{2}\langle X(y-x), y-x\rangle
       +o(|s-t|+|y-x|^2).
$$ Similarly, we denote by $\Cal P^{2,-}u(t,x)$ [the "parabolic
subjet" of $u$ at $(t,x)$] the set of triples
$(p,q,X)\in\IR\times\IR^n\times S^n$ which are such that
$$
u(s,y) \geq u(t,x)+p(s-t)+\langle q, y-x\rangle
        +\frac{1}{2}\langle X(y-x), y-x\rangle
       +o(|s-t|+|y-x|^2).
$$}

\no{\bf Example 4.2}$\quad$ {\it  Suppose that $\varphi\in
C^{1,2}((0,T)\times\IR^n)$. If $u-\varphi$ has a local maximum at
$(t,x)$, then
$$\left(\frac{\partial \varphi}{\partial t}(t,x),\nabla\varphi(t,x),D^2\varphi(t,x)\right)\in\Cal P^{2,+}u(t,x).$$
If $u-\varphi$ has a local minimum at $(t,x)$, then
$$\left(\frac{\partial \varphi}{\partial t}(t,x),\nabla\varphi(t,x),D^2\varphi(t,x)\right)\in\Cal P^{2,-}u(t,x).$$
}\med

We can now give the definition of a viscosity solution of the HJB
equation obstacle problem (4.1) .\med

\no{\bf Definition 4.3}

(a) {\it It can be said $u(t,x)\in C([0,T]\times \IR^n)$ is a
viscosity subsolution of (4.1) if $u(T,x)\leq \Phi(x)$,
$x\in\IR^n$, and at any point $(t,x)\in (0,T)\times\IR^n$, for any
$(p,q,X)\in\Cal P^{2,+}u(t,x)$,
$$\min\left(u(t,x)-h(t,x),-p-\sup_{v\in U}\left\{\frac 1 2 Tr(aX)+
\langle
b,q\rangle+g(t,x,u(t,x),q\sigma(t,x,v),v)\right\}\right)\leq 0.$$
In other words at any point $(t,x)$ where $u(t,x)>h(t,x)$},
$$-p-\sup_{v\in U}\left\{\frac 1 2 Tr(aX)+\langle b,q\rangle+g(t,x,u(t,x),q\sigma(t,x,v),v)\right\}\leq 0.$$

 (b) {\it  It can be said $u(t,x)\in C([0,T]\times \IR^n)$ is a
viscosity supersolution of (4.1) if $u(T,x)\geq \Phi(x)$,
$x\in\IR^n$, and at any point $(t,x)\in (0,T)\times\IR^n$, for any
$(p,q,X)\in\Cal P^{2,-}u(t,x)$,
$$\min\left(u(t,x)-h(t,x),-p-\sup_{v\in U}\left\{\frac 1 2 Tr(aX)+\langle b,q\rangle+g(t,x,u(t,x),q\sigma(t,x,v),v)\right\}\right)\geq 0.$$
In other words, at each point, we have both $u(t,x)\geq h(t,x)$
and}
$$-p-\sup_{v\in U}\left\{\frac 1 2 Tr(aX)+\langle b,q\rangle+g(t,x,u(t,x),q\sigma(t,x,v),v)\right\}\geq 0.$$

 (c) {\it  $u(t,x)\in C([0,T]\times \IR^n)$ is said to be a viscosity
solution of (4.1) if it is both a viscosity sub- and
supersolution. } \med

We are going to use the approximation of the reflected BSDE by
penalization, which was studied in section 6 of [9]. For each
$(t,x)\in [0,T]\times\IR^n$, $n\in\IN$, let
$\{(^nY^{t,x;v}_s,^nZ^{t,x;v}_s), t\leq s\leq T\}$ denote the
solution of the BSDE
$$\aligned
^nY^{t,x;v}_s &= \Phi(X^{t,x;v}_T)+\int_s^T
g(r,X^{t,x;v}_r,^nY^{t,x;v}_r,^nZ^{t,x;v}_r,v_r)dr\\
&\quad +n\int_s^T(^nY^{t,x;v}_r-h(r,X^{t,x;v}_r))^-dr-\int_s^T\
^nZ^{t,x;v}_rdW_r,\quad t\leq s\leq T.
\endaligned$$
We define
$$
J_n(t,x;v(\cdot)) := ^nY^{t,x;v}_t,\qquad v(\cdot)\in\Cal U,\
0\leq t\leq T,\ x\in\IR^n;\tag 4.2$$
$$ u_n(t,x) :=
ess\sup_{v(\cdot)\in\Cal U} J_n(t,x;v(\cdot)),\qquad 0\leq t\leq
T,\ x\in\IR^n.\tag 4.3
$$
It is known from [12] or [13]  that $u_n(t,x)$ defined in (4.3) is
the viscosity solution of the PDE
$$\left\{\aligned &-\frac{\partial u_n}{\partial t}(t,x)-\sup_{v\in U}\left\{\Cal L(t,x,v)u_n(t,x)
+g_n(t,x,u_n(t,x),\nabla u_n(t,x)\sigma(t,x,v),v)\right\}=0,\\
&u_n(T,x)=\Phi(x),\endaligned\right.$$ where
$$g_n(t,x,r,p\sigma(t,x,v),v)=g(t,x,r,p\sigma(t,x,v),v)+n(r-h(t,x))^-.$$
Then

\no{\bf Lemma 4.4}   $u_n(t,x)\uparrow u(t,x)$, $0\leq t\leq T $,
$x\in\IR^n$. \med

{\bf Proof}:   From the result of the section 6 in [9], for each
$0\leq t\leq T$, $x\in\IR^n$,
$$J_n(t,x;v(\cdot))\uparrow J(t,x;v(\cdot)),\qquad {\text as }\quad  n\to\infty.$$
From the monotonic property of $J_n$ and the definition of $u_n$
in (4.3), we get the monotonic property of $u_n$. Next we will
show the convergent property of $u_n$.

For each $0\leq t\leq T$, $x\in\IR^n$, $\forall\varepsilon>0$,
there exists $v(\cdot)\in\Cal U$ such that
$$u(t,x)<Y^{t,x;v}_t+\varepsilon,$$
then
$$0\leq u(t,x)-u_n(t,x)\leq Y^{t,x;v}_t- ^nY^{t,x;v}_t+\varepsilon.$$
Because $^nY^{t,x;v}_t\uparrow Y^{t,x;v}_t$, $a.s.$, we take limit
on both side,
$$0\leq \limsup_{n\to\infty}(u(t,x)-u_n(t,x))\leq \varepsilon.$$
From the arbitrariness of $\varepsilon$, we get the desired
result.
\endpf\med

\no{\bf Remark 4.5}  {\it Since $u_n$ and $u$ are continuous, it
follows from Dini's theorem that the convergence in the lemma is
uniform on compacts.}\med

\no{\bf Theorem 4.6}  {\it Defined by (3.9), $u$ is a viscosity
solution of HJB equations (4.1).}

\med {\bf Proof}:  We now show that $u$ is a subsolution of (4.1).
Let $(t,x)$ be a point at which $u(t,x)>h(t,x)$, and let
$(p,q,X)\in\Cal P^{2,+}u(t,x)$.

From Lemma 6.1 in [4], there exists sequences
$$n_j\to+\infty,\quad
(t_j,x_j)\to (t,x),\quad (p_j,q_j,X_j)\in\Cal
P^{2,+}u_{n_j}(t_j,x_j),$$ such that
$$(p_j,q_j,X_j)\to (p,q,X).$$
But for any $j$,
$$\aligned
& -p_j-\sup_{v\in U}\left\{\frac{1}{2}Tr(aX_j)+\langle
b,q_j\rangle+g(t_j,x_j,u_{n_j}(t_j,x_j),q_j\sigma(t_j,x_j,v),v)\right.\\
&\qquad\qquad
\left.+n_j(u_{n_j}(t_j,x_j)-h(t_j,x_j))^-\right\}\leq 0.
\endaligned$$
From the assumption that $u(t,x)>h(t,x)$ and the uniform
convergence of $u_n$, it follows that for $j$ large enough
$u_{n_j}(t_j,x_j)>h(t_j,x_j)$, hence
$$-p_j-\sup_{v\in U}\left\{\frac{1}{2}Tr(aX_j)+\langle
b,q_j\rangle+g(t_j,x_j,u_{n_j}(t_j,x_j),q_j\sigma(t_j,x_j,v),v)\right\}\leq
0.$$

Let us admit for a moment the following lemma.\med

\no{\bf Lemma 4.7 }
$$\aligned
& \lim_{j\to\infty}\sup_{v\in U}\left\{\frac{1}{2}Tr(aX_j)+\langle b, q_j\rangle+g(t_j,x_j,u_{n_j}(t_j,x_j),q_j\sigma(t_j,x_j,v),v)\right\}\\
&= \sup_{v\in U}\lim_{j\to\infty}\left\{\frac{1}{2}Tr(aX_j)+\langle
b,
q_j\rangle+g(t_j,x_j,u_{n_j}(t_j,x_j),q_j\sigma(t_j,x_j,v),v)\right\}.
\endaligned$$

Taking the limit as $j\to\infty$ in the above inequality yields:
$$-p-\sup_{v\in U}\left\{\frac{1}{2}Tr(aX)+\langle
b,q\rangle+g(t,x,u(t,x),q\sigma(t,x,v),v)\right\}\leq 0,$$ and we
have proved that $u$ is a subsolution of (4.1).\med

We now show that $u$ is a supersolution of (4.1). Let $(t,x)$ be
an arbitrary point in $(0,T)\times\IR^n$, and $(p,q,X)\in\Cal
P^{2,-}u(t,x)$. We already know that $u(t,x)\geq h(t,x)$. By the
same argument as above, there exist sequences:
$$n_j\to+\infty,\quad (t_j,x_j)\to (t,x),\quad
(p_j,q_j,X_j)\in\Cal P^{2,-}u_{n_j}(t_j,x_j),$$ such that
$$(p_j,q_j,X_j)\to (p,q,X).$$
But for any $j$,
$$\aligned
& -p_j-\sup_{v\in U}\left\{\frac{1}{2}Tr(aX_j)+\langle
b,q_j\rangle+g(t_j,x_j,u_{n_j}(t_j,x_j),q_j\sigma(t_j,x_j,v),v)\right.\\
&\qquad\qquad
\left.+n_j(u_{n_j}(t_j,x_j)-h(t_j,x_j))^-\right\}\geq 0.
\endaligned$$
Hence
$$-p_j-\sup_{v\in U}\left\{\frac{1}{2}Tr(aX_j)+\langle
b,q_j\rangle+g(t_j,x_j,u_{n_j}(t_j,x_j),q_j\sigma(t_j,x_j,v),v)\right\}\geq
0,$$ and taking the limit as $j\to\infty$, we conclude that:
$$-p-\sup_{v\in U}\left\{\frac{1}{2}Tr(aX)+\langle
b,q\rangle+g(t,x,u(t,x),q\sigma(t,,v),v)\right\}\geq 0.$$
\endpf

Now we turn to \med

\no{\bf  Proof of Lemma 4.7} $\quad$ For the convenience, we
denote
$$f_j(v)=\frac{1}{2}Tr(a,X_j)+\langle b,q_j\rangle+g(t_j,x_j,u_{n_j}(t_j,x_j),q_j\sigma(t_j,x_j,v),v).$$
Firstly, $\forall v\in U$,
$$f_j(v)\leq\sup_{v\in U}f_j(v),\quad
\lim_{j\to\infty}f_j(v)\leq \liminf_{j\to\infty}\sup_{v\in
U}f_j(v),$$ then
$$\sup_{v\in U}\lim_{j\to\infty}f_j(v)\leq \liminf_{j\to\infty}\sup_{v\in U}f_j(v).\tag 4.4$$

Secondly, we consider a subsequence $\{j_k\}_{k=1}^\infty$ such
that
$$\lim_{j_k\to\infty}\sup_{v\in U}f_{j_k}(v)=\limsup_{j\to\infty}\sup_{v\in U}f_j(v).$$
$\forall\varepsilon>0$, $\forall j_k$, $\exists v_{j_k}\in U$ such
that
$$\sup_{v\in U}f_{j_k}(v)\leq f_{j_k}(v_{j_k})+\varepsilon.$$
Because $U$ is compact, there exists a convergent subsequence
denoted by $\{v_{j_k}\}_{k=1}^\infty$ also, the limit is denoted
by $v_0$. We consider the difference of $f_{j_k}(v_{j_k})$ and
$f_{j_k}(v_0)$: From the Lipschitz condition we get
$$|f_{j_k}(v_{j_k})-f_{j_k}(v_0)|\leq C|v_{j_k}-v_0|^2+C|v_{j_k}-v_0|,$$
where $C$ only depend on the Lipschitz constant. It follows that
for $j_k$ large enough
$$|f_{j_k}(v_{j_k})-f_{j_k}(v_0)|\leq \varepsilon.$$
Then
$$\sup_{v\in U}f_{j_k}(v)\leq f_{j_k}(v_0)+2\varepsilon,$$
$$\limsup_{j\to\infty}\sup_{v\in U}f_j(v)=\lim_{j_k\to\infty}\sup_{v\in U}f_{j_k}(v)\leq\lim_{j_k\to\infty}f_{j_k}(v_0)+2\varepsilon=\lim_{j\to\infty}f_j(v_0)+2\varepsilon,$$
$$\limsup_{j\to\infty}\sup_{v\in U}f_j(v)\leq\sup_{v\in U}\lim_{j\to\infty}f_j(v_0)+2\varepsilon.$$
From the arbitrariness of $\varepsilon$,
$$\limsup_{j\to\infty}\sup_{v\in U}f_j(v)\leq\sup_{v\in U}\lim_{j\to\infty}f_j(v_0).\tag 4.5$$
From (4.4) and (4.5), we complete the proof. \endpf
\med

Finally, we shall use some technique and method from [1] to
establish a uniqueness result for viscosity solution of (4.1). This
kind of technique and method can also be seen in [3] to prove the
uniqueness for viscosity solutions of Hamilton-Jacobi-Bellman-Isaacs
equations related to stochastic differential games.\med

\no{\bf Lemma 4.8} $\quad$  {\it Let $u_1\in C([0,T]\times\IR^n)$
be a viscosity subsolution and $u_2\in C([0,T]\times\IR^n)$ be a
viscosity supersolution of (4.1). Then the function $w:=u_1-u_2$
is a viscosity subsolution of the system
$$\left\{\aligned &\min\left(w(t,x),-\frac{\partial w}{\partial t}(t,x)-\sup_{v\in U}\left\{\Cal L(t,x,v)w(t,x)+L|w|
+L|\nabla w\sigma(t,x,v)|\right\}\right)=0,\\
&w(T,x)=0,\endaligned\right.\tag 4.6$$ where $L$ is the Lipschitz
constant of $g$ in $(y,z)$. } \med

{\bf Proof}: The proof is similar to that of the corresponding
results: Lemma 3.7 in [1].\med

For each $(t_0,x_0)\in (0,T)\times\IR^n$, let $\varphi\in
C^{\infty}([0,T]\times\IR^n)$ and let $(t_0,x_0)$ be a strict
global maximum point of $w-\varphi$. Because $u_2$ is a viscosity
supersolution of HJB equation (4.1), we have $u_2(t_0,x_0)\geq
h(t_0,x_0)$. If $u_1(t_0,x_0)\leq h(t_0,x_0)$, it is easily to get
$$w(t_0,x_0)=u_1(t_0,x_0)-u_2(t_0,x_0)\leq 0,$$
and we get the desired result. Therefore, in the proof, we suppose
that $u(t_0,x_0)>h(t_0,x_0)$.

We introduce the function
$$\Phi_\varepsilon(t,x,y)=u_1(t,x)-u_2(t,y)-\frac{|x-y|^2}{\varepsilon^2}-\varphi(t,x),$$
where $\varepsilon$ is a positive parameter which is devoted to
tend to zero.

Since $(t_0,x_0)$ is a strict global maximum point of
$u_1-u_2-\varphi$, by a classical argument in the theory of
viscosity solutions, there exists a sequence
$(\hat{t},\hat{x},\hat{y})$ such that\med

(i) $(\hat{t},\hat{x},\hat{y})$ is a global maximum point of
$\Phi_\varepsilon$ in $[0,T]\times\bar{B}_R\times\bar{B}_R$ where
$B_R$ is a ball with a large radius $R$;

(ii) $(\hat{t},\hat{x})$, $(\hat{t},\hat{y})\to (t_0,x_0)$ as
$\varepsilon\to 0^+$;

(iii) $\frac{|\hat{x}-\hat{y}|^2}{\varepsilon^2}$ is bounded and
tend to zero when $\varepsilon\to 0^+$.\med

We have dropped above the dependence of $\hat{t}$, $\hat{x}$ and
$\hat{y}$ in $\varepsilon$ for the sake of simplicity of
notations.

It follows from Theorem 8.3 in [4] that, $\forall\delta>0$, there
exist
$$
\left(p,\frac{2(\hat{x}-\hat{y})}{\varepsilon^2}+D\varphi,X\right)\in\bar{\Cal
P}^{2,+}u_1(\hat{t},\hat{x}),\quad \left(p-\frac{\partial
\varphi}{\partial
t},\frac{2(\hat{x}-\hat{y})}{\varepsilon^2},Y\right)\in \bar{\Cal
P}^{2,-}u_2(\hat{t},\hat{y}),$$ such that
$$\left(\matrix  X & 0 \\ 0 & -Y \endmatrix\right)\leq A+\delta A^2,\tag 4.7$$
where
$$A=\left(\matrix \frac{2}{\varepsilon^2}+D^2\varphi & -\frac{2}{\varepsilon^2}
 \\ -\frac{2}{\varepsilon^2} & \frac{2}{\varepsilon^2} \endmatrix\right).$$
Calculating directly, we get
$$\aligned
 A+\delta A^2 & =
\left(\frac{2}{\varepsilon^2}+\delta\frac{4}{\varepsilon^4}\right)\left(\matrix
I & -I \\ -I & I \endmatrix\right)
+\left(1+\delta\frac{4}{\varepsilon^2}\right)\left(\matrix
D^2\varphi & 0 \\ 0 & 0 \endmatrix\right)\\
&\quad +\delta\frac{4}{\varepsilon^4}\left(\matrix I & 0 \\
0 & I \endmatrix\right) +\delta\left(\matrix (D^2\varphi)^2 & 0
\\ 0 & 0 \endmatrix\right).
\endaligned$$
After given $\varepsilon$, $\delta>0$, we have
$$\aligned
&
-p-\sup_{v\in U}\left\{\frac{1}{2}Tr\left((\sigma\sigma^T)(\hat{t},\hat{x},v)X\right)+
\langle b(\hat{t},\hat{x},v),\frac{2(\hat{x}-\hat{y})}{\varepsilon^2}+D\varphi(\hat{t},\hat{x})\rangle\right.\\
&\hskip 2cm \left.+g\left(\hat{t},\hat{x},u_1(\hat{t},\hat{x}),[\frac{2(\hat{x}-\hat{y})}{\varepsilon^2}+D\varphi(\hat{t},\hat{x})]\sigma(\hat{t},\hat{x},v),v\right)\right\}\leq 0,\\
&
-\left(p-\frac{\partial \varphi}{\partial t}(\hat{t},\hat{x})\right)-\sup_{v\in U}\left
\{\frac{1}{2}Tr\left((\sigma\sigma^T)(\hat{t},\hat{y},v)Y\right)
+\langle b(\hat{t},\hat{y},v),\frac{2(\hat{x}-\hat{y})}{\varepsilon^2}\rangle\right.\\
&\hskip 2cm
\left.+g\left(\hat{t},\hat{y},u_2(\hat{t},\hat{y}),\frac{2(\hat{x}-\hat{y})}{\varepsilon^2}\sigma(\hat{t},\hat{y},v),v\right)\right\}\geq
0.
\endaligned$$
The first inequality minus the second one,
$$\aligned
& -\frac{\partial \varphi}{\partial t}(\hat{t},\hat{x})-\sup_{v\in U}\left
\{\frac{1}{2}\left(Tr\left((\sigma\sigma^T)(\hat{t},\hat{x},v)X\right)
-Tr\left((\sigma\sigma^T)(\hat{t},\hat{y},v)Y\right)\right)\right.\\
&\quad +\left(\langle
b(\hat{t},\hat{x},v),\frac{2(\hat{x}-\hat{y})}{\varepsilon^2}
+D\varphi(\hat{t},\hat{x})\rangle-\langle b(\hat{t},\hat{y},v),\frac{2(\hat{x}-\hat{y})}{\varepsilon^2}\rangle\right)\\
&\quad
\left.+\left[g\left(\hat{t},\hat{x},u_1(\hat{t},\hat{x}),[\frac{2(\hat{x}-\hat{y})}{\varepsilon^2}+D\varphi(\hat{t},\hat{x})]\sigma(\hat{t},\hat{x},v)\right)
         -g\left(\hat{t},\hat{y},u_2(\hat{t},\hat{y}),\frac{2(\hat{x}-\hat{y})}{\varepsilon^2}\sigma(\hat{t},\hat{y},v)\right)\right]\right\}\\
&\leq 0.
\endaligned$$
Using (4.7) and Lipschitz condition, we analysis the items in the
$\sup_{v\in U}$ and  get
$$\aligned
&  -\frac{\partial \varphi}{\partial t}(\hat{t},\hat{x})
       -\sup_{v\in U}\left\{\frac{1}{2}(\frac{2}{\varepsilon^2}+\delta\frac{4}{\varepsilon^4})L^2|\hat{x}
       -\hat{y}|^2+\frac{1}{2}(1+\delta\frac{4}{\varepsilon^2})Tr\left((\sigma\sigma^T)(\hat{t},\hat{x},v)D^2
       \varphi(\hat{t},\hat{x})\right)\right.\\
&\quad
+\frac{1}{2}\delta\frac{4}{\varepsilon^4}(|\sigma(\hat{t},\hat{x},v)|^2+|\sigma(\hat{t},\hat{y},v)|^2)
                           +\frac{1}{2}\delta Tr\left((\sigma\sigma^T)(\hat{t},\hat{x},v)(D^2\varphi)^2(\hat{t},\hat{x})\right)\\
&\quad  + 2L \frac{|\hat{x}-\hat{y}|^2}{\varepsilon^2}   +\langle
b(\hat{t},\hat{x},v),D\varphi(\hat{t},\hat{x})\rangle
+L|\hat{x}-\hat{y}|+L|u_1(\hat{t},\hat{x})-u_2(\hat{t},\hat{x})|\\
&\quad \left. +L|u_2(\hat{t},\hat{x})-u_2(\hat{t},\hat{y})|
+L|D\varphi(\hat{t},\hat{x})\sigma(\hat{t},\hat{x},v)|+2L^2\frac{|\hat{x}-\hat{y}|^2}{\varepsilon^2}\right\}
\leq 0.
\endaligned$$
We let $\delta\to 0^+$, then  let $\varepsilon\to 0^+$ and  we get
$$\aligned
 -\frac{\partial \varphi}{\partial t}(t_0,x_0)
                                  &-\sup_{v\in U}\left\{\frac{1}{2}Tr\left((\sigma\sigma^T)(t_0,x_0,v)D^2
                                  \varphi(t_0,x_0)\right)\right.+\langle  b(t_0,x_0,v),D\varphi(t_0,x_0)\rangle\\
&\quad  \left.
+L|w(t_0,x_0)|+L|D\varphi(t_0,x_0)\sigma(t_0,x_0,v)|\right\}\leq
0.
\endaligned$$
Therefore $w$ is a viscosity subsolution of the desired equation
(4.6) and  the proof is completed.
\endpf
\med

Now we are going to construct one suitable smooth supersolution
for the equation (4.6).\med

\no{\bf Lemma 4.9}  {\it For any $A>0$, there exists $C_1>0$ such
that the function
$$\chi(t,x)=\exp\left\{(C_1(T-t)+A)\psi(x)\right\},$$
where
$$\psi(x)=\left[\log\left((|x|^2+1)^{\frac{1}{2}}\right)+1\right]^2$$
satisfies
$$\min\left(\chi(t,x),-\frac{\partial \chi}{\partial t}(t,x)-\sup_{v\in U}\left\{\Cal L(t,x,v)\chi(t,x)+L\chi(t,x)+L|\nabla \chi\sigma(t,x,v)|\right\}\right)>0$$
in $[t_1,T]\times\IR^n$ where $t_1=T-(A/C_1)$.
 }
\med

{\bf Proof}:  Obviously, the function $\chi$ defined in the Lemma
satisfy $\chi(t,x)>0$, for each $(t,x)\in [0,T]\times \IR^n$. We
give estimations on the first and second order derivatives of
$\psi$:
$$|D\psi(x)|\leq\frac{2[\psi(x)]^{\frac{1}{2}}}{(|x|^2+1)^{\frac{1}{2}}}\quad
{\text and} \quad
|D^2\psi(x)|\leq\frac{C\left(1+[\psi(x)]^{\frac{1}{2}}\right)}{|x|^2+1}\quad
{\text in}\quad \IR^n.$$ These estimations imply that, if $t\in
[t_1,T]$,
$$|D\chi(t,x)|\leq C\chi(t,x)\frac{[\psi(x)]^{\frac{1}{2}}}{(|x|^2+1)^{\frac{1}{2}}},\quad
|D^2\chi(t,x)|\leq C\chi(t,x)\frac{\psi(x)}{|x|^2+1},$$
where the constant $C$ only depend on $A$. We continue to
calculate
$$\aligned
&\frac{\partial \chi}{\partial t}(t,x)+\sup_{v\in U}\left\{\Cal
L(t,x,v)\chi(t,x)+L\chi(t,x)+L|\nabla
\chi\sigma(t,x,v)|\right\}\\
&=\frac{\partial \chi}{\partial t}(t,x)+\sup_{v\in
U}\left\{\frac{1}{2}Tr((\sigma\sigma^T)D^2\chi)+\langle b,
D\chi\rangle+L\chi(t,x)+L|\nabla
\chi\sigma(t,x,v)|\right\}\\
&\leq -C_1\chi(t,x)\psi(x)+\sup_{v\in U}\left\{\frac 1 2
\frac{|\sigma(t,x,v)|^2}{|x|^2+1}C\chi(t,x)\psi(x)\right.\\
&\quad\left.+\frac{|b(t,x,v)|}{(|x|^2+1)^{\frac{1}{2}}}C\chi(t,x)[\psi(x)]^{\frac{1}{2}}
+L\chi(t,x)+L\frac{|\sigma(t,x,v)|}{(|x|^2+1)^{\frac{1}{2}}}C\chi(t,x)[\psi(x)]^{\frac{1}{2}}\right\}.
\endaligned\tag 4.8$$
Because $b$ and $\sigma$ are linear growth in $x$,
$[\psi(x)]^{\frac{1}{2}}\leq\psi(x)$ and $1\leq\psi(x)$, the above
inequality (4.8)
$$\aligned
&< -C_1\chi(t,x)\psi(x)+\frac 1 2
C\chi(t,x)\psi(x)+C\chi(t,x)\psi(x)+L\chi(t,x)\psi(x)+LC
\chi(t,x)\psi(x)\\
&=-(C_1-\frac 1 2 C-C-L-LC)\chi(t,x)\psi(x).
\endaligned$$
It is clear that when $C_1$ large enough the quantity in the right
side of the above inequality is negative and the proof is
completed.
\endpf
\med

Now we can prove the uniqueness result for viscosity solution of
(4.1).\med

\no{\bf Theorem 4.10}  {\it Assume that $b$, $\sigma$, $g$, $\Phi$
and $h$ satisfy (H3.1)--(H3.4), respectively. Then there exists at
most one viscosity solution of HJB equation (4.1) in the class of
continuous functions which grow at most polynomially at infinity.}
\med

{\bf Proof}:  Let $u_1, u_2\in C([0,T]\times\IR^n)$ be two
viscosity solutions of HJB equation (4.1).

We define $w:=u_1-u_2$, then we have
$$\lim_{|x|\to\infty}w(t,x)e^{-A[\log((|x|^2+1)^{\frac{1}{2}})]^2}=0$$
uniformly for $t\in [0,T]$, for some $A>0$. This implies, in
particular, that $w(t,x)-\alpha\chi(t,x)$ is bounded from above in
$[t_1,T]\times\IR^n$ for any $\alpha>0$ and that
$$M:=\max_{[t_1,T]\times\IR^n}(w-\alpha\chi)(t,x)e^{-L(T-t)}$$
is achieved at some point $(t_0,x_0)\in [t_1,T]\times\IR^n$
(depend on $\alpha$). Then we have two case.\med

The first case: $w(t_0,x_0)\leq 0$.

Then we have
$$u_1(t,x)-u_2(t,x)\leq \alpha\chi(t,x),\qquad (t,x)\in [t_1,T]\times\IR^n.$$
Letting $\alpha$ tends to zero, we obtain
$$u_1(t,x)\leq u_2(t,x),\qquad (t,x)\in [t_1,T]\times\IR^n.\tag 4.9$$
\med

The second case: $w(t_0,x_0)>0$.

Then we have
$$w(t,x)-\alpha\chi(t,x)\leq (w(t_0,x_0)-\alpha\chi(t_0,x_0))e^{-L(t-t_0)},\qquad (t,x)\in [t_1,T]\times\IR^n.$$
We define
$$\varphi(t,x)=\alpha\chi(t,x)+(w(t_0,x_0)-\alpha\chi(t_0,x_0))e^{-L(t-t_0)},$$
and can get
$$w-\varphi\leq 0=(w-\varphi)(t_0,x_0),\qquad (t,x)\in [t_1,T]\times \IR^n.$$
Since $\varphi(t_0,x_0)=w(t_0,x_0)>0$ and Lemma 4.8, when $t_0\in
[t_1,T)$, we have
$$\aligned
&-\frac{\partial \varphi}{\partial t}(t_0,x_0)-\sup_{v\in
U}\left\{\frac{1}{2}Tr\left((\sigma\sigma^T)(t_0,x_0,v)D^2\varphi(t_0,x_0)\right)+\langle
b(t_0,x_0,v), D\varphi(t_0,x_0)\rangle\right.\\
&\quad
+L\varphi(t_0,x_0)+L|\nabla\varphi(t_0,x_0)\sigma(t_0,x_0,v)|\big\}\leq
0.
\endaligned$$
From the definition of $\varphi$, we rewrite the above inequality
$$\aligned
&\alpha\left [-\frac{\partial \chi}{\partial
t}(t_0,x_0)-\sup_{v\in
U}\left\{\frac{1}{2}Tr\left((\sigma\sigma^T)(t_0,x_0,v)D^2\chi(t_0,x_0)\right)+\langle
b(t_0,x_0,v), D\chi(t_0,x_0)\rangle\right.\right.\\
&\quad\left.+L\chi(t_0,x_0)+L|\nabla\chi(t_0,x_0)\sigma(t_0,x_0,v)|\big\}\right]\leq
0.
\endaligned$$
This is a contradiction with Lemma 4.9. Therefore $t_0=T$, this is
a contradiction with the fact that $w(t,x)$ is a viscosity
subsolution of (4.6) (see Lemma 4.8). Then the second case does
not happen.\med

If we change $w(t,x)=u_1-u_2$ for $w'(t,x)=u_2-u_1$, the same
argument leads to
$$u_2(t,x)\leq u_1(t,x),\qquad (t,x)\in [t_1,T]\times\IR^n.\tag 4.10$$
Combining  (4.9) with (4.10), we have
$$u_1(t,x)=u_2(t,x),\qquad (t,x)\in [t_1,T]\times\IR^n.$$
\med
Applying successively the same argument on the intervals
$[t_2,t_1]$ where $t_2=(t_1-A/C_1)^+$ and then, if $t_2>0$ on
$[t_3,t_2]$ where $t_3=(t_2-A/C_1)^+$ ... etc. We finally obtain
that
$$u_1(t,x)=u_2(t,x),\qquad (t,x)\in [0,T]\times\IR^n.$$
The proof is complete.
\endpf

\med \centerline{\bf  Appendix} \med

In the appendix we give the proof of Proposition 2.1 and 2.2. \med
\no{\bf Proof of Proposition 2.1} \med

Applying It\^{o}'s formula to the process $|Y_s|^2e^{\beta s}$
yields
$$\aligned
&|Y_t|^2e^{\beta t}+\int_t^T(\beta|Y_s|^2+|Z_s|^2)e^{\beta
s}ds \\
&= |\xi|^2e^{\beta T}+2\int_t^T Y_sg(s,Y_s,Z_s)e^{\beta s}ds
+2\int_t^TY_se^{\beta s}dK_s-2\int_t^T Y_sZ_se^{\beta
s}dW_s\\
&=|\xi|^2e^{\beta T}+2\int_t^T Y_sg(s,Y_s,Z_s)e^{\beta s}ds
+2\int_t^TS_se^{\beta s}dK_s-2\int_t^T Y_sZ_se^{\beta s}dW_s,
\endaligned$$
where we have used the identity $\int_t^T(Y_s-S_s)e^{\beta
s}dK_s=0$. Using the Lipschitz property of $g$, we have
$$\aligned
&|Y_t|^2e^{\beta t}+\int_t^T(\beta|Y_s|^2+|Z_s|^2)e^{\beta
s}ds \\
&\leq |\xi|^2e^{\beta T}+2\int_t^T |Y_s||g(s,Y_s,Z_s)|e^{\beta s}ds
+2\int_t^T |S_s|e^{\beta s}dK_s-2\int_t^T Y_sZ_se^{\beta
s}dW_s\\
&\leq |\xi|^2e^{\beta T}+2\int_t^T |Y_s||g(s,0,0)|e^{\beta s}ds
+2\int_t^T (L|Y_s|^2+L|Y_s||Z_s|)e^{\beta s}ds\\
&\quad +2\int_t^T |S_s|e^{\beta s}dK_s-2\int_t^T Y_sZ_se^{\beta
s}dW_s\\
&\leq |\xi|^2e^{\beta T}+2\int_t^T |Y_s||g(s,0,0)|e^{\beta s}ds
+\int_t^T \left((2L+2L^2)|Y_s|^2+\frac{1}{2}|Z_s|^2\right)e^{\beta s}ds\\
&\quad +2\int_t^T |S_s|e^{\beta s}dK_s-2\int_t^T Y_sZ_se^{\beta
s}dW_s.
\endaligned$$
We select $\beta=2L^2+2L$, then
$$\aligned
|Y_t|^2e^{\beta t}+\frac 1 2 \int_t^T|Z_s|^2e^{\beta s}ds &\leq
|\xi|^2e^{\beta T}+2\int_t^T |Y_s||g(s,0,0)|e^{\beta s}ds\\
&+2\int_t^T |S_s|e^{\beta s}dK_s-2\int_t^T Y_sZ_se^{\beta
s}dW_s.\endaligned\tag A.1$$
$$\aligned
\IE^{\Cal F_t}\left\{\int_t^T|Z_s|^2e^{\beta s}ds\right\}&\leq
2\IE^{\Cal F_t}\left\{|\xi|^2e^{\beta T}+2\int_t^T
|Y_s||g(s,0,0)|e^{\beta s}ds\right.\\
&\left.+2\int_t^T |S_s|e^{\beta s}dK_s\right\}.\endaligned\tag A.2$$
$$\aligned \sup_{t\leq u\leq T} |Y_u|^2e^{\beta u} &\leq
|\xi|^2e^{\beta T}
+2\int_t^T|Y_s||g(s,0,0)|e^{\beta s}ds\\
&\quad +2\int_t^T |S_s|e^{\beta s}dK_s +4\sup_{t\leq u\leq
T}|\int_t^uY_sZ_se^{\beta s}dW_s|.
\endaligned$$
From Burkholder-Davis-Gundy's inequality we have
$$\aligned
\IE^{\Cal F_t}\left\{\sup_{t\leq u\leq T} |Y_u|^2e^{\beta u}\right\}
&\leq \IE^{\Cal F_t}\left\{|\xi|^2e^{\beta T}
+2\int_t^T|Y_s||g(s,0,0)|e^{\beta s}ds\right.\\
&\left.+2\int_t^T |S_s|e^{\beta s}dK_s\right\}+C\IE^{\Cal
F_t}\left(\int_t^T|Y_s|^2|Z_s|^2e^{2\beta s}ds\right)^{\frac 1 2},
\endaligned$$
 thanks to the inequality $ab\leq a^2/2+b^2/2$, we deduce
immediately $$\aligned \IE^{\Cal F_t}\left\{\sup_{t\leq u\leq T}
|Y_u|^2e^{\beta u}\right\} &\leq \IE^{\Cal
F_t}\left\{|\xi|^2e^{\beta T} +2\int_t^T|Y_s||g(s,0,0)|e^{\beta
s}ds+2\int_t^T |S_s|e^{\beta s}dK_s \right\}\\
&\quad+\frac{C^2}{2}\IE^{\Cal F_t}\left\{\int_t^T|Z_s|^2e^{\beta
s}ds\right\}+\frac{1}{2}\IE^{\Cal F_t}\left\{\sup_{t\leq u\leq
T}|Y_u|^2e^{\beta u}\right\}.
\endaligned$$
Combining the inequality (A.2) with the above one, we easily derive
that
$$\aligned
 &\IE^{\Cal F_t}\left\{\sup_{t\leq u\leq
T}|Y_u|^2e^{\beta u}
+\int_t^T|Z_s|^2e^{\beta s}ds\right\}\\
&\leq C\IE^{\Cal F_t}\left\{|\xi|^2e^{\beta T}
+\int_t^T|Y_s||g(s,0,0)|e^{\beta s}ds +2\int_t^T |S_s|e^{\beta
s}dK_s \right\}.
\endaligned$$
Using the fact that $$ C\IE^{\Cal
F_t}\left\{\int_t^T|Y_s||g(s,0,0)|e^{\beta s}ds\right\} \leq
\frac{1}{2}\IE^{\Cal F_t}\left\{\sup_{t\leq u\leq T}|Y_u|^2e^{\beta
u}\right\}  +\frac{C^2}{2}\IE^{\Cal
F_t}\left(\int_t^T|g(s,0,0)|e^{(\beta/2)s}ds\right)^2, $$ we get
$$\aligned &\IE^{\Cal F_t}\left\{\sup_{t\leq u\leq
T}|Y_u|^2e^{\beta u}
+\int_t^T|Z_s|^2e^{\beta s}ds\right\}\\
&\leq C\IE^{\Cal F_t}\left\{|\xi|^2e^{\beta T}
+\left(\int_t^T|g(s,0,0)|e^{(\beta/2)s}ds\right)^2 +2\int_t^T
|S_s|e^{\beta s}dK_s \right\}.
\endaligned$$
Then we drop the exponential function to get a brief form
$$\aligned
&\IE^{\Cal F_t}\left\{\sup_{t\leq u\leq T}|Y_u|^2
+\int_t^T|Z_s|^2ds\right\}\\
&\quad \leq C\IE^{\Cal F_t}\left\{|\xi|^2
+\left(\int_t^T|g(s,0,0)|ds\right)^2 +2\int_t^T |S_s|dK_s
\right\}.\endaligned\tag A.3
$$
We now give an estimate of $\IE^{\Cal F_t}[|K_T-K_t|^2]$. From the
equation
$$K_T-K_t=Y_t-\xi-\int_t^T g(s,Y_s,Z_s)ds+\int_t^T Z_sdW_s,$$
and estimate (A.3), we get the following inequalities:
$$\aligned
 \IE^{\Cal F_t}\left\{|K_T-K_t|^2\right\} &\leq C\IE^{\Cal
F_t}\left\{|\xi|^2 +\left(\int_t^T|g(s,0,0)|ds\right)^2 +2\int_t^T
|S_s|dK_s\right\}\\
&\leq C\IE^{\Cal F_t}\left\{|\xi|^2
+\left(\int_t^T|g(s,0,0)|ds\right)^2 \right\}\\
&\quad +2C^2\IE^{\Cal F_t}\left\{\sup_{t\leq u\leq T}
|S_s|^2\right\} +\frac{1}{2}\IE^{\Cal
F_t}\left\{|K_T-K_t|^2\right\}.
\endaligned$$
Consequently,
$$\IE^{\Cal F_t}\left\{|K_T-K_t|^2\right\}\leq C\IE^{\Cal F_t}\left\{|\xi|^2
+\left(\int_t^T|g(s,0,0)|ds\right)^2 +\sup_{t\leq u\leq T}
|S_s|^2\right\}.\tag A.4$$
 Combining the estimate (A.3) with
(A.4), we complete the proof of the proposition.
\endpf
\med

\no{\bf Proof of Proposition 2.2} \med

 The computation process is
similar to that in the  proof of Proposition 2.1, so we shall only
give the sketch of the proof. Since $\int_t^T(\Delta Y_s-\Delta
S_s)e^{\beta s}d(\Delta K_s)\leq 0$,
$$\aligned
& |\Delta Y_t|e^{\beta t}+\int_t^T(\beta|\Delta Y_s|^2+|\Delta
Z_s|^2)e^{\beta s}ds\\
&\leq |\Delta \xi|^2e^{\beta T}+2\int_t^T\Delta Y_s\Delta
g(s,Y_s,Z_s)e^{\beta s}ds\\
&\quad +2\int_t^T\Delta Y_s[g'(s,Y_s,Z_s)-g'(s,Y'_s,Z'_s)]e^{\beta
s}ds\\
&\quad +2\int_t^T\Delta S_se^{\beta s}d(\Delta K_s)-2\int_t^T\Delta
Y_s \Delta Z_se^{\beta s}dW_s\\
&\leq |\Delta \xi|^2e^{\beta T}+2\int_t^T|\Delta Y_s||\Delta
g(s,Y_s,Z_s)|e^{\beta s}ds\\
&\quad +2L\int_t^T(|\Delta Y_s|^2+|\Delta Y_s||\Delta Z_s|)e^{\beta
s}ds\\
&\quad +2\int_t^T|\Delta S_s|e^{\beta s}d(K_s+K'_s)-2\int_t^T\Delta
Y_s \Delta Z_se^{\beta s}dW_s.
\endaligned$$
Similar technique with the above proof of Proposition 2.1, we can
get $$\aligned &\IE^{\Cal F_t}\left\{\sup_{t\leq u\leq T}|\Delta
Y_u|^2+\int_t^T|\Delta Z_s|^2ds\right\}\\
&\leq C\IE^{\Cal F_t}\left\{|\Delta \xi|^2+\left(\int_t^T|\Delta g(s,Y_s,Z_s)|ds\right)^2+2\int_t^T|\Delta S_s|d(K_s+K'_s)\right\}\\
&\leq C\IE^{\Cal F_t}\left\{|\Delta \xi|^2+\left(\int_t^T|\Delta
g(s,Y_s,Z_s)|ds\right)^2+\left(\sup_{t\leq u\leq T}|\Delta S_u|\right)\left((K_T-K_t)+(K'_T-K'_t)\right)\right\}\\
&\leq C\IE^{\Cal F_t}\left\{|\Delta \xi|^2+\left(\int_t^T|\Delta g(s,Y_s,Z_s)|ds\right)^2\right\}\\
&\quad +\left(\IE^{\Cal F_t}\left\{\sup_{t\leq u\leq T} |\Delta
S_u|^2\right\}\right)^{1/2} \left(\IE^{\Cal
F_t}\left\{((K_T-K_t)-(K'_T-K'_t))^2\right\}\right)^{1/2}.
\endaligned$$
And then from Proposition 2.1, we complete the proof.
\endpf

\med \no{\bf Acknowledgements.} The authors express  their gratitude
to Prof. Shige Peng for his elicitation and inspiring idea in
recursive stochastic dynamic programming principle. The authors also
thank Dr. Juan Li and Mingyu Xu for their helpful discussions and
suggestions.

 \med \no{\bf Reference} \med \item{[1]}G. Barles, R. Buckdahn \& E. Pardoux,
Backward Stochastic Differential Equations and Integral-Partial
Differential Equations, Stochastics and Stochastics Reports,
60(1997), pp. 57-83.

\med \item{[2] } P. Briand, F. Coquet, Y. Hu, J. M\'{e}min \& S.
Peng, A converse comparison theorem for BSDEs and related properties
of $g-$expectation, Elect. Comm. in Probab., 5(2000), pp. 101-117.

\med \item{[3]} R. Buckdahn \& J. Li, Stochastic Differential games
and Viscosity solutions of Hamilton-Jacobi-Bellman-Isaacs Equations,
Preprint, 2006.

\med\item{[4]} M. G. Crandall, H. Ishii \& P. L. Lions,  User's
guide to viscosity solutions of second order partial differential
equations,
 Bull. Amer. Soc.,  27(1992), pp. 1-67.
 \med
\item{[5] } J.Cvitanic \& I.Karatzas,  Backward SDE's with
reflection and Dynkin Games, The Annals of Probability,  24(1996),
pp. 2024-2056. \med \item{[6]} D.~Duffie \& L.~Epstein,  Stochastic
differential utility, Econometrica,  60(1992), pp. 353--394. \med
\item{[7] } S.Hamad\`ene \& J.-P.Lepeltier,  Reflected BSDEs and
mixed  game problems, Stochastic processes and their applications, {
85(2000), pp. 177-188. \med \item{[8] } S.Hamad\`ene, J.-P.Lepeltier
\& Z.Wu,  Infinite horizon Reflected BSDEs and applications in mixed
control and game problems, Probability and mathematical statistics,
 19(1999), pp. 211-234.

 \med\item{[9] } El.
~Karoui, C.~Kapoudjian, E.~Pardoux, S.~Peng \& M.C.~Quenez ,
Reflected solutions of backward SDE's, and related obstacle problems
for PDE's, The Annals of Probability,  25(1997), pp. 702-737. \med
\item{[10]}N.~El.Karoui, S.~Peng \& M.C.~Quenez,  Backward
Stochastic Differential Equation in Finance, Math. Finance, 7(1997),
pp. 1-71. \med \item{[11] } E.Pardoux \& S.Peng,  Adapted solutions
of a backward
 stochastic differential equation,
Systems and Control Letters,  14(1990), pp. 55-61. \med
\item{[12] } S.Peng,  A generalized dynamic programming
principle and Hamilton-Jacobi-Bellmen equation, Stochastics and
Stochastic Reports,  38(1992), pp. 119-134. \med\item{[13]} J. Yan,
S.Peng, S.Fang \& L.Wu,  Topics on stochastic analysis, Science
Press. Beijing (in Chinese), 1997. \med

\enddocument